\documentclass{article}
\usepackage{amsmath}
\usepackage{amssymb}
\usepackage{graphicx}
\usepackage{multirow}
\usepackage{algorithm}
\usepackage{algorithmic}
\usepackage{float}
\usepackage{url}
\usepackage{lipsum}
\usepackage{multicol}

\usepackage{tikz}
\usetikzlibrary{shapes.geometric, arrows}
\tikzstyle{startstop} = [rectangle, rounded corners, minimum width=3cm, minimum height=1cm,text centered, draw=black, fill=red!30]
\tikzstyle{io} = [trapezium, trapezium left angle=70, trapezium right angle=110, minimum width=3cm, minimum height=1cm, text centered, draw=black, fill=blue!30]
\tikzstyle{process} = [rectangle, minimum width=3cm, minimum height=1cm, text centered, draw=black, fill=orange!30]
\tikzstyle{decision} = [diamond, minimum width=3cm, minimum height=1cm, text centered, draw=black, fill=green!30]
\tikzstyle{arrow} = [thick,->,>=stealth]

\newtheorem{thm}{\bf Theorem}[section]

\newtheorem{exmp}[thm]{\bf Example}

\title{A Simplified Algorithm for Identifying Abnormal Changes in Dynamic Networks}

\date{}
\author{}

\begin{document}

 \thispagestyle{empty}
%%%%%%%%%%%%%%%%%%%%%%%%%%%%%%%%%%%%%%%%%%%%%%%%%%%%%%%%%
%%%%%%%%%%%%%%%%%%%%%%%%%%%%%%%%%%%%%%%%%%%%%%%%%%%%%%%%%

%%%%%%%%%%%%%%%%%%%%%%%%%%%%%%%%%%%%%%%%%%%%%%%%%%%%%%%%%
%%%TITLE%%%%%%%%%%%%%%%%%%%%%%%%%%%%%%%%%%%%%%%%%%%%%%%%%
\maketitle \vspace*{-1.5cm}

  \begin{center}
{\large\bf Bouchaib Azamir$^{1, a}$,   Driss Bennis$^{1, b}$ and Bertrand Michel$^{2, c}$}

\bigskip

%%%%%%%%%%%%%%%%%%%%%%%%%%%%%%%%%%%%%%%%%%%%%%%%%%%%%%%%%
%%%ADDRESSES%%%%%%%%%%%%%%%%%%%%%%%%%%%%%%%%%%%%%%%%%%%%%

$^1$   Faculty of Sciences,  Mohammed V University in Rabat,  Morocco.\\
\noindent    $^a$\,bouchaib\_azamir@um5.ac.ma \\ $^b$\,driss.bennis@um5.ac.ma; driss$\_$bennis@hotmail.com  \\[0.2cm]

$^2$   Nantes Université, Ecole Centrale Nantes, Laboratoire de mathématiques Jean Leray UMR 6629.\\ 1 Rue de La Noe, 44300 Nantes, France.\\
\noindent  $^c$\,bertrand.michel@ec-nantes.fr\bigskip
\end{center}

%
%%%%%%%
%%%ABSTRACT%%%%%%%%%%%%%%%%%%%%%%%%%%%%%%%%%%%%%%%%%%%%%%
\noindent{\large\bf Abstract.} 
 Topological data analysis has recently been applied to the study of dynamic networks. In this context, an algorithm was  introduced and helps, among other things, to detect early warning signals of abnormal changes in the dynamic network under study. However, the complexity of this algorithm increases significantly once the database studied grows. In this paper, we propose a simplification of the algorithm without affecting its performance.  We give various applications and simulations of the new algorithm on some weighted networks. The obtained results   show clearly the efficiency of the introduced approach. Moreover, in some cases, the proposed algorithm makes it possible to highlight local information and sometimes early warning signals of local abnormal changes.
\bigskip


\section{Introduction}

A network is an efficient representation of a given set of entities and the existing relationships between them. Actually a network is a graph which entities are the vertices and the relations between them are presented as edges (see for instance  \cite{19} and \cite{4}).  Usually, these relations are quantified by (mostly positive) real numbers and in this case we say that the graph (the network) is weighted.
In most real-life situations, (weighted) networks evolve over time resulting in a family of (weighted) networks; we then speak of dynamic (weighted) networks (see the beginning of Section 3). In general, there are two kinds of dynamic networks : 
\begin{itemize}
\item There are dynamic networks such that both the set of vertices and the set of edges are changing over time. For instance, chat groups in social networks where the members of these groups are undergoing changes from moment to moment.
\item In some situations, the set of vertices remains unchanged over time as in the case of financial networks such as the weighted networks representing stock market transactions. So, the only changes occur at the level of the edges and at the level of the weights.
\end{itemize}

Dynamic networks have been subject of several studies. See for instance the papers and the references within \cite{Li4,Li5,Li3,Li1,Li2} where dynamic networks have been studied following different approaches. In this paper, we are interested in studying the behaviour of dynamic weighted networks using   topological data analysis (TDA). In \cite{11}, Gidea used TDA to detect the early signs of a financial market crisis. His method is based mainly on the fact that we can associate with any financial network of stocks a metric space whose distance is defined by the correlation coefficients of stock market returns. This allows to compute the so-called persistent homology  of this metric space (see the preliminaries in Section 2)  and, thanks to the concept of time series, the first signs of a critical transition in the financial network are detected. It is clear that Gidea's method can be applied on any other dynamic network. However, for large networks, the execution time of the algorithms used in Gidea's method will be very high. Thus, one can ask whether the fact of working on particular subgraphs rather than the whole network permits to achieve good results in a reasonable amount of time. Indeed, the abundance of real-world data and its velocity usually call for data summarization. Namely, graph summarization is a useful tool to help identify the structure and meaning of data. It has various advantages including reducing data volume and storage, speeding up algorithms and graph queries, eliminating noise, among others. There are different types of graph summaries depending on whether the data is homogeneous or heterogeneous and also whether the network is static or dynamic. In the literature, graph summarization methods use a set of basic techniques, here we are interested in those that are based on simplification or parsimony: these methods rationalize an input graph by removing nodes or edges less ``important", resulting in a sparse graph see \cite{21}.  Here, the summary graph consists of a subset of the original nodes and/or edges. A representative work on node simplification-based summarization techniques is OntoVis \cite{22}, a visual analysis tool for exploring and understanding large heterogeneous social networks. In \cite{23}, the authors propose a four-step unsupervised algorithm for egocentric information abstraction of heterogeneous social networks using edge, instead of node, filtering. In  \cite{24}, a new problem, to simplify weighted graphs by removing less important edges, has been proposed. Simplified graphs can be used to improve the visualization of a network, to extract its main structure, or as a preprocessing step for other data mining algorithms. 
Our approach is part of the simplification of graphs and aims to identify a subnetwork that contains a large part of the information embodied by the original network.
\\

Our aim in this paper is to show that there are some subnetworks that can be good candidates for this study. Those subnetworks are extracted from the whole graph by eliminating some edges at certain threshold  which is determined by the use of the closeness centrality concept (see Section \ref{Sec-Central} for the definition of closeness centrality). 
The idea is based on the fact that  a node having the highest closeness centrality is  able to spread information efficiently through the whole graph (see \cite{9} for more details). Then, based on this fact, we propose that  the desired subnetwork would be the one which contains this node called henceforth central node. After that, we need to set the threshold where edges with higher weights will be eliminated. Also, following the closeness centrality notion,  the threshold will be chosen from the list of weights of the incident edges to the central node. However, it is not clear how to determine the optimal threshold. This is why we propose to focus our study on the three quartiles, the minimum and the maximum of the above proposed list of weights. The associated  subnetworks will be called central subnetworks (see  Section \ref{Sec-Central} for more details about centrality).  Our aim is to determine what threshold levels provides good results. We will see throught several examples that when the threshold is greater than  the third quartiles, the adopted approach reduces considerably the execution time of Gidea's method (to almost half in some cases)  without compromising quality and results.  Moreover, in some cases we get good results even with thresholds greater than the median.
 \\

This paper is organized as follows:\\

Section \ref{prelimin} presents some preliminaries on persistent homology. We start this  section with a part (Subsection \ref{sub-sec-simplicial}) in which we present a short  mathematical background that we need to understand the persistent homology associated to   weighted graphs (given in Subsection \ref{PH}), and we  present also how to implement  data to get   persistence diagrams associated to weighted graphs (see Subsection \ref{implem}). 

 In Section \ref{Sec3}, we  present how to associate persistent homology to dynamic networks (see Subsection \ref{description}) and show how it was used by Gidea to detect abnormal changes in financial networks  (see Subsection \ref{sec-Gidea-appli}).  

Section  \ref{Sec-Central} represents our approach to reduce the computation cost of Gidea's method (see Figure \ref{flowchart} for the 
 flowchart of the proposed method). It is mainly based on a simplification of the whole studied network by considering what we called central subnetworks determined using the notion of closeness centrality (see Algorithm 1). We also describe how to get the persistence diagrams associated to these central subnetworks (see Algorithm 2). 
  
Section   \ref{sec-experimnts} concerns the application phase. Namely, we present some numerical experiments in which we compare our method with existing ones. In each of the first three subsections, we simulate a dynamic network.   Each subsection deals with a kind of a  weighted network whose vertex point cloud is simulated with a given probability distribution. Then the corresponding central subnetworks are derived. A visual comparison between the time series associated to both the whole network and the central subnetwork shows that are very similar (see Figures in this section). This observation is also confirmed by the calculated R-squared coefficients  (see Tables in this section). In parallel,  we show that using our method the execution time is reduced considerably (see Tables  in this section).  
We end Section   \ref{sec-experimnts} by a comparison a  performance comparison between our approach and the edge collapse method (Subsection \ref{compa-edge}).  It is known that  the edge collapse method provides theoretical guarantees and allows to obtain, in output, the same starting persistence diagram (see for instance \cite{collapse}). From Table \ref{tbl:compa-e-collaps}, we deduce that our method is more practical for large and dense networks. 

We end this article with Section \ref{appli dynamic} which includes an application of the method on real weighted networks. In Subsection \ref{appli DJIA}, we apply on  the financial network used by Gidea \cite{11},  and in Subsection \ref{crypto} we consider a cryptocurrency network. As in Section   \ref{sec-experimnts}, we deduce an improvement of the execution time using the new method.

%%%%%%%%%%%%%%%%%%%%%%%%%%%%%%%%%%%%%%%%%%%%%%%%%%%%%%%%%%%%%%%%%%%%%%%%%%%%%%%%%%%%%%%%%%%%%%%%%%%%%%%%%%%%%%%%%%%%%%%%%%%%%%%%%%%%%%%%%%%%%%%%%%%%%%

%%%%%%%%%%%%%%%%%%%%%%%%%%%%%%%%%%%%%%%%%%%%%%%%%%%%%%%%%%%%%%%%%%%%%%%%%%%%%%%%%%%%%%%%%%%%%%%%%%%%%%%%%%%%%%%%%%%%%%%%%%%%%%%%%%%%%%%%%%%%%%%%%%%%%%
\section{Persistent homology computed on graphs }\label{prelimin}
Persistent homology is an emergent method for detecting geometrical and topological properties of a space endowed with a topological structure. The concept of persistence was first introduced by Edelsbrunner, Letscher and Zomorodian in \cite{8} and then refined by Carlsson and  Zomorodian in  \cite{2}. Since then, persistent homology has become an essential tool in topological Data Analysis (TDA) and has been applied in various scientific fields such as biology, image processing, sensor networks
etc. Nowadays, there are many references  (surveys, books and notes) on TDA; for a recent one we refer  to  \cite{5}.\\

 In this section, we give a short presentation on some basic notions concerning particular simplicial complexes as well as its corresponding  persistent homology. We will focus, following our context, on simplicial complexes associated to graphs investigated in \cite{1}. For a general background on persistent homology  one can see  \cite{2, 5, 7,10}. 

 \subsection{Simplicial complex and   simplicial homology groups}\label{sub-sec-simplicial}
We start this subsection with a recall of a simplicial complex and then we give the examples of interest.\\

Given a finite set $ V$, a simplicial complex with the vertex set $V$ is a set $\widetilde{K}$ of finite subsets of $V$ such that the elements of $V$ belong to $\widetilde{K}$ and for any $\sigma \in \widetilde{K}$, any subset of $\sigma$ belongs to $\widetilde{K}$. The elements of $\widetilde{K}$ are called the faces or the simplices of $\widetilde{K}$. The dimension of a simplex of   $\widetilde{K}$   is just its cardinality minus 1 and the dimension of the simplicial complex $\widetilde{K}$ is the largest dimension of its simplices \cite{5}.

There are a lot of examples of simplicial complexes arising from various contexts. Here, we use the following ones:

\begin{exmp}[Simplicial clique complex,  \cite{1}]
Let $G=(V,E)$ be a graph, where $V$ is the set of vertices and $E$ the set of edges. For a positive integer $k$, a  $k$-clique of $G$ is a set of $k$ vertices  whose the induced subgraph is complete. Denote by $Cl (G)$ the set of all the cliques of $G$. Notice that the set  $Cl (G)$ satisfies the  two conditions :
\begin{itemize}
  \item $Cl (G)$ contains all singletons $\{v\}$ with $v\in V$.
  \item  $Cl (G)$  is closed under subsets: if $\tau \subseteq \sigma  $ and $\sigma \in Cl (G) $, then $\tau \in Cl (G)$.
\end{itemize}
Then, $Cl (G)$ is  a simplicial complex which is called a simplicial clique complex of $G$.
\end{exmp}

The following (geometric) classical example of simplicial complexes  can be found in any introductory books on TDA.

\begin{exmp}[Vietoris-Rips complex]
Given a finite set of points $\mathbb{X}$ in a metric space $ (M, d)$ and a real number $\alpha \ge 0$. The Vietoris-Rips complex $\mathrm{Rips}_{\alpha}(X)$ is a simplicial complex whose simplices are sets   $\{x_ 0, … , x_ k \}$ such that $d(x_i,x_j) \le \alpha $ for all $ 0\leq i, j\leq k$.
\end{exmp}

Now we are ready to recall what are simplicial homology groups of  a simplicial complex. To this end, we define some vector spaces and linear maps.
Here, we  only deal with vector spaces on the field   $\mathbb{F}_2:=\mathbb{Z}/2\mathbb{Z}$. \\

Given a simplicial complex  $K$, the  vector space generated by the $p$-simplices of $K$ is denoted by $Cp(K)$. It consists of all finite formal sums of $p$-simplices called $p$-chains; i.e., an element $c$ belongs to $Cp(K)$ if it can be written as follows: $ c=\sum_{j}\gamma_{j}\sigma_{j} $ for some scalers $\gamma_{j}\in \mathbb{F}_2$ and  a family $(\sigma_{j})_j$  of $p$-simplices.

For a positive integer $p$, consider the linear map $\partial_p: C_{p}(K) \rightarrow C_{p-1}(K)$, called boundary map, defined on $p$-simplices as follows: For every $p$-simplex $\sigma$, $\partial(\sigma)$ is the formal sum of the $(p-1)$-dimensional faces (i.e., subsets of $\sigma$ of cardinal $p$). An element of the image of $\partial_p$ is called a boundary. Thus, the boundary $\partial_p(c)$ of a chain $ c=\sum_{j}\gamma_{j}\sigma_{j} $   is obtained by extending $\partial_p$ linearly; i.e.,  $\partial_p(c)=\sum_{j}\gamma_{j}\partial_p(\sigma_{j}) $.
   The $p$-chains that have boundary $0$ are called $p$-cycles. They form a subspace $Z_{p}$ of $C_{p}$. The $p$-chains which are the boundary of $(p + 1)$-chains are called $p$-boundaries and form a subspace $B_{p}$ of $C_{p}$. It is important but not difficult to show that $\partial_p\circ\partial_{p+1}= 0 $ which is equivalent to  $B_{p}\subseteq Z_{p}$. Then,  if we consider  the quotient  vector space  $Z_{p}/B_{p}$, that means we have annihilated the   $p$-boundaries which are not  $p$-cycles. For instance, if we take $K$ to be $Cl (G)$, the simplicial clique complex of a graph $G$, the dimension of the vector space $B_{1}/ Z_{1}$ will be  exactly the number of ``holes"  in $G$ and the dimension of the vector space $B_{0}/ Z_{0}$ is the number of the connected components of  $Cl (G)$ (i.e., the connected components of  the graph (network) $G$).  For this reason and of course other important ones,     the vector space   $B_{p}/ Z_{p}$ is considered as one of the principal notions  in algebraic topology  that helps   in  describing topological features. It is denoted by $ H_{p}(K)$   and usually called the  $p^{th}$ simplicial homology group of $K$. Its elements are called homology classes.

  \begin{exmp}\label{exm-Cl}
   The simplicial clique complex $Cl (G)$ associated to the graph $G$ in Figure \ref{fig-graph_example} has three connected components, and so the dimension of $ H_{0}(Cl (G))$ is three. But it has only one hole because $A_{9}A_{10}+$ $A_{10}A_{11}+$ $A_{11}A_{12}+A_{12}A_{9}$ is a cycle but not a boundary. Notice that $A_{5}A_{6}+A_{6}A_{7}+A_{7}A_{8}+A_{8}A_{5}$ is a boundary of the $3$-chain $A_{5}A_{6}A_{8}+A_{6}A_{7}A_{8}$.
\begin{figure}[H]
\centering
  \fbox{  \includegraphics[scale=0.8]{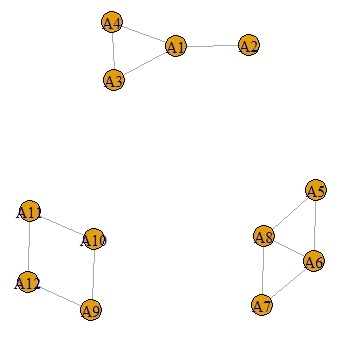}}
  \caption{A graph $G$ with three connected components.}\label{fig-graph_example}
\end{figure} 
  \end{exmp}

%%%%%%%%%%%%%%%%%%%%%%%%%%%%%%%%%%%%%%%%%%%%%%%%
%%%%%%%%%%%%%%%%%%%%%%%%%%%%%%%%%%%%%%%%%%%%%%%%%%%%%%%%%%%%%%%%%%%%%%%%

\subsection{Persistence diagrams}\label{PH}
 When a simplicial complex $K$  has a history in the sense that it is an union of a chain of  simplicial complexes, we can capture more information on $K$; especially, the behaviour of homology groups of these complexes can be used as a signature of the complex $K$ that can be used to distinguish it from  other complexes.  This is one of many important things that persistent homology  is performing. \\

\begin{figure}[H]
  \centering
 \fbox{  \includegraphics[scale=0.4]{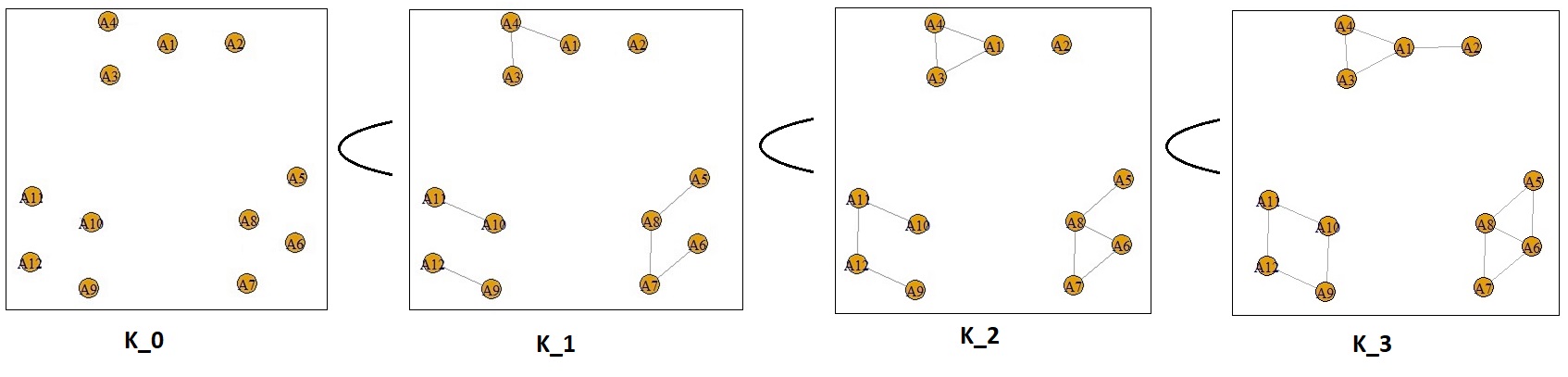}}
  \caption{Filtration of the simplicial clique complex   given in Example \ref{exm-Cl}.}\label{fig-filtration_example}
\end{figure}

Consider the following  chain of   simplicial subcomplexes of a complex $  K$, called a filtration of $K$:
$$\emptyset\subseteq K_0\subseteq K_1\subseteq \cdots \subseteq K_p=K  $$
  This filtration provides more information on K in the following sense:  for every couple   $i,j$ such that $ 0\leq i\leq j\leq p$, the injection  $K_i \hookrightarrow K_j$ induces a homomorphism on the  $n^{th}$ simplicial homology groups for each dimension:
	$$f_{n}^{i,j}:H_n(K_i)\rightarrow H_n(K_j).$$
    The $n^{th} $ persistent Betti number $\beta_{n}^{i,j}$ is the rank of the vector space $Im(f_{n}^{i,j})$; that is,  $\beta_{n}^{i,j} = dim(Im(f_{n}^{i,j}))$. Persistent Betti numbers count how many homology classes of dimension $n$ that survive during the passage from $K_i$ to $K_j$. We say that a homology class $\alpha \in H_n(K_i)$ is born entering in $K_i$, if $\alpha$ does not come from a previous subcomplex; that is, $\alpha \notin Im(f_{n}^{i-1,i})$. Similarly, if $\alpha$ is born in $K_i$, it dies entering $ K_j$ if the image of the map induced by $ K_{i-1}\subseteq K_{j-1}$ does not contain the image of $\alpha$  but the image of the map induced by $ K_{i-1}\subseteq K_j $ does. In this case the persistence of $\alpha$ is $j-i$. Finally, all the collected data of births and deaths are presented in   a diagram  which is  the multi-set of points with  coordinates   correspond respectively to the birth time and the death time of  homology classes. This diagram is called an $n$-dimensional persistence diagram, or simply a persistence diagram, of $K$.\\
As a simple example, the filtration represented in Figure \ref{fig-filtration_example} is a filtration of the simplicial clique complex  $Cl (G)$  given in Example \ref{exm-Cl}.

\begin{figure}[H]
  \centering
 \fbox{  \includegraphics[scale=0.4]{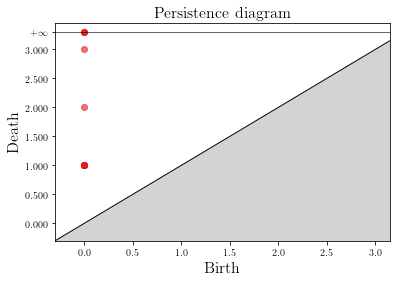}}
  \caption{The ($0$-dimensional) persistence diagram associated to the filtration of Figure  \ref{fig-filtration_example}.}\label{fig-diag_example}
\end{figure}

The $0$-dimensional persistence diagram associated to the filtration given in  Figure  \ref{fig-filtration_example} is given in Figure \ref{fig-diag_example}.

Certainly, the persistence diagram of $ K $ reveals more information on it than those extracted from the (static) complex $ K $.  Moreover, persistence diagrams can be used to compute the degree of similarity between complexes. To this end,   some metric distances are introduced on  the space of persistence diagrams.  In this paper, we use the so called $p$-Wasserstein distance (for a positive integer $p $). It is defined by the formula:   $$ W_p (\Delta_1, \Delta_2 )=\underset{\phi}{\inf} [\sum_{q\in \Delta_1}\|q-\phi(q)\|]^{1/p}  $$
where      $\Delta_1$ and $ \Delta_2$ are two diagrams and $\phi$ varies across all one-to-one and onto functions from $\Delta_1$ to $ \Delta_2$. In this paper we use just the $2$-Wasserstein distance.\\

Now, to apply the concept above on graph,   the  notion of 
 weight functions is used as follows:  Let $G=(V,E,\omega)$ be a weighted graphs, where $V$ is the set of vertices and $E$ the set of edges,    such that the  weight function  $\omega:E\rightarrow \mathbb{R}^+$ is assumed to be upper bounded by a positive real number $m$. Let $K=Cl(G)$ be the clique complex of $G$. We would like to construct a filtration of $K$ (i.e., a chain of complexes that covers $K$). To do so, let us consider,  for every  $a\in \mathbb{R}^+$,   the sublevel set $ E(a)=\omega^{-1} ((-\infty,a])$ and the associated subgraph $G(a)$ of $G$   which will be   called the subgraph of $G$ at the  threshold $a$.  Then, we consider  the associated clique complex  $ K(a)$ at the $a^{th} $ level. In particular, $ K(m)=K$. Now, for a sequence $ a_0<a_1<\cdots<a_p=m  $  of real numbers, we get a filtration  of $K$ (called Rips filtration):
$$\emptyset\subseteq K_0\subseteq K_1\subseteq \cdots \subseteq K_p=K  $$
where $ K_i=K(a_i)$.  Therefore,    the persistent homology, described above, applied on  this filtration provides more information on the graph $G$ as we will deduce throughout this paper.

%%%%%%%%%%%%%%%%%%%%%%%%%%%%%%%%%%%%%%%%%%%%%%%%
%%%%%%%%%%%%%%%%%%%%%%%%%%%%%%%%%%%%%%%%%%%%%%%%%%%%%%%%%%%%%%%%%%%%%%%%

\subsection{Implementation}\label{implem}
Nowadays, there are various  softwares  for   computations of persistent  homology and corresponding information. In this article, we often  use the R package TDA  (for a short tutorial and introduction on the use of the R package TDA, see \cite{FKL}), except in Subsection \ref{compa-edge} where we use also the Gudhi library\footnote{https://gudhi.inria.fr}. We will  also use the R software for computations corresponding the   time series (see for instance \cite{18}).

In this subsection, we present  how to implement the method  described in Subsection  \ref{PH} in order to get the persistence diagram of a weighted graph.
It should be noted that the existing software  implements data from a metric space. Thus, in order to implement the method described in Subsection \ref{PH}, we need a distance between nodes.  The weight function can be used to this end (see Subsection \ref{sec-Gidea-appli}). So here we assume that the  weight function    defines a distance between nodes. In other words, the graph with the weight  function is a metric space. In this case, the associated clique complexes are nothing but the  Vietoris-Rips complexes.\\

Now, the implementation which provides persistence diagrams of a  weighted graph, as it is conceived,  deals with matrices of distances associated to a metric space. These matrices   are, in our context,   matrices of weighted complete graphs. 
Then, given a weighted complete graph $G$, one can determine its persistence diagram by using the Rips filtration method such that the max-scale of the filtration is the maximum of the weights.

It is also possible to compute the persistence diagram of a not necessary complete graph. Let $G = (V,E,\omega)$ be a weighted graph, where $V$ and $E$ are the sets of vertices and edges, respectively,  and $\omega:E\to \mathbb{R^+}$ is a weight function  assumed to be upper bounded by a positive real number $m$. The data of this graph can be represented in a square   matrix $W=(W_{i,j})$ such that $ W_{i,j} =\omega(i,j)$ if $(i,j)\in E$ and $ W_{i,j} = m+1$ otherwise. This approach is adopted to not take in consideration the non connected vertices in the  filtration process  since the corresponding sequence of real numbers, as explained in Subsection \ref{PH}, ends at $m$.\\
Notice that if one would like to determine the persistence diagram of $G(s)$, the  subgraph of  $G$ at a threshold $s$  strictly less than $m$,
 we consider the matrix $\widetilde{W}=(\widetilde{W}_{i,j})$  such that  $ \widetilde{W}_{i,j} =W_{i,j}$ if $W_{i,j}\leq s$ and $ \widetilde{W}_{i,j} = m+1$ otherwise.

%%%%%%%%%%%%%%%%%%%%%%%%%%%%%%%%%%%%%%%%%%%%%%%%%%%%  To construct a filtration of $K=Cl(G)$, we consider

%%%%%%%%%%%%%%%%%%%%%%%%%%%%%%%%%%%%%%%%%%%%%%%%%%%%

\section{Persistent homology of dynamic networks}\label{Sec3}

In this section, we present how to associate persistent homology to dynamic networks and show how it was used by Gidea to detect abnormal changes in financial networks.

   \subsection{Description of the  method}\label{description}
    First, we recall what is a dynamic network.
As concluded in \cite[Section 2.4]{4},  there is no consensual formal definition of a dynamic network. Some papers define a dynamic system as a distributed system in which the communication graph evolves over time, while others define it as a model where processes can join and leave the system during a run. Some articles mix both of those proposals in the same model (\cite{4}).  As   presented in  \cite{19}: a (discrete) dynamic weighted network can be mathematically  represented as a time sequence of weighted graphs,  $(G_t = (V_ t,E_t,\omega_t))_{t\in \mathbb{N}}$, where $V_t$ and $E_t$ are the set of vertices and the set of edges, respectively, and $\omega_t:E_t\to \mathbb{R^+}$ is a weight function at time $t$. If the set of vertices remains unchanged over time; that is, $V_ t=V$ for all $t$,  then the weighted graph at any time $t$ will be simply   $G_t = (V, E_t,\omega_t)$. In this case, a dynamic weighted network  can be encoded as a  weighted  matrix, $W(t)=(W_{i,j}(t))$, where $W_{i,j}(t)$ is the  weight of the edge between vertices $i$ and $j$ at time  $t$.\\

Now, by applying the method described in Subsection \ref{PH} on each weighted graph $G_t$, we end up with a times series  of persistence diagrams $( \Delta_t )_t$.  To induce a scalar time series from $( \Delta_t )_t$, we consider a  fixed   persistence diagram $\Delta$,  then the    $2$-Wasserstein distances,
$W_2 (\Delta, \Delta_t )$,
between the persistence diagrams $ \Delta_t $ and $\Delta$, give rise to a (scalar) time series.
 This time series  makes it possible to detect abnormal changes within the dynamic network.   This idea was used by  Gidea in \cite{11} to  detect abnormal changes  in a financial network.    The following subsection presents briefly Gidea's work.

    \subsection{Gidea's application}\label{sec-Gidea-appli}
    The dynamic weighted network studied by Gidea in \cite{11}  was a financial network with a weight function defined via correlations between nodes (stocks) as follows:
     For each stock $i$ and a day $t$, $x_i (t)$ is the daily returns based on the adjusted closing prices $ S_i (t)$;  i.e., $x_i (t)=\frac{S_i (t+1)-S_i (t)}{S_i (t)}$.
  Let $ x_i (t) $ be the arithmetic return of a stock $ i $ at time $t$. The correlation coefficient between nodes $i$ and $j$ over a time interval of the form $ [t-T, T]$, where $T>0$,     is defined by:	$$ C_{i,j}(t):=\frac{\sum_{\tau ={t-T}}^{T}(x_{i}(\tau)-\overline{x_i})(x_{j}(\tau)-\overline{x_j})}{\sqrt{\sum_{\tau=t-T}^{T}(x_i(\tau)-\overline{x_i})^2}{{\sqrt{\sum_{\tau=t-T}^{T}(x_j(\tau)-\overline{x_j})^2}}}}$$
       where $\overline{x_i}$ and
   $\overline{x_j} $ denote the averages of $x_j $ and $ x_j $, respectively, over the time interval $ [t-T,T] $.   Then, following \cite[Section 13.1]{14}, $ d(i,j)(t)= \sqrt{2(1-C_{i,j}(t))} $ is a distance  function  between nodes $i$ and $j$. Notice that the values taken by the distance $d$ vary between 0 and 2; especially,   if two nodes $i$ and $j$ are perfectly correlated we have $d(i,j)(t)=0$, while $d(i,j)(t)=2$  when they are perfectly anti-correlated. With  $d$,  this network is a metric space and so it can be seen as a time evolving weighted network; i.e.,  as  a graph $ G(V, E)$ with the set $V$ of nodes representing the collection of stocks and  the weight function  $ \omega_t$ (at a time $t$) from the set of edges $E$ into $ [0, \infty)$ defined by $ \omega_t (e)=d(i,j)(t)$ for every edge $e$. Thus, we can generate the associated persistence diagram $\Delta_t$ (at a time $t$) as described in Section \ref{prelimin}.   We know   that the set of these  persistence diagrams is  endowed with the $2$-Wasserstein distance. We fix a persistence diagram $\Delta_{t_0}$ at an initial time $t_0$, then the behaviour of the time series  $(X_t)_t:=( W_2 (\Delta_t, \Delta_{t_0} ) )_t$  will translate the (topological) changes in the studied  financial network. In particular,   it could help  in detecting a change prior to the critical transition (i.e., the peak of the crisis),  meaning that the stock correlation network undergoes significant changes in its topological structure.\\

\begin{figure}[H]
  \centering
  \includegraphics[scale=0.07]{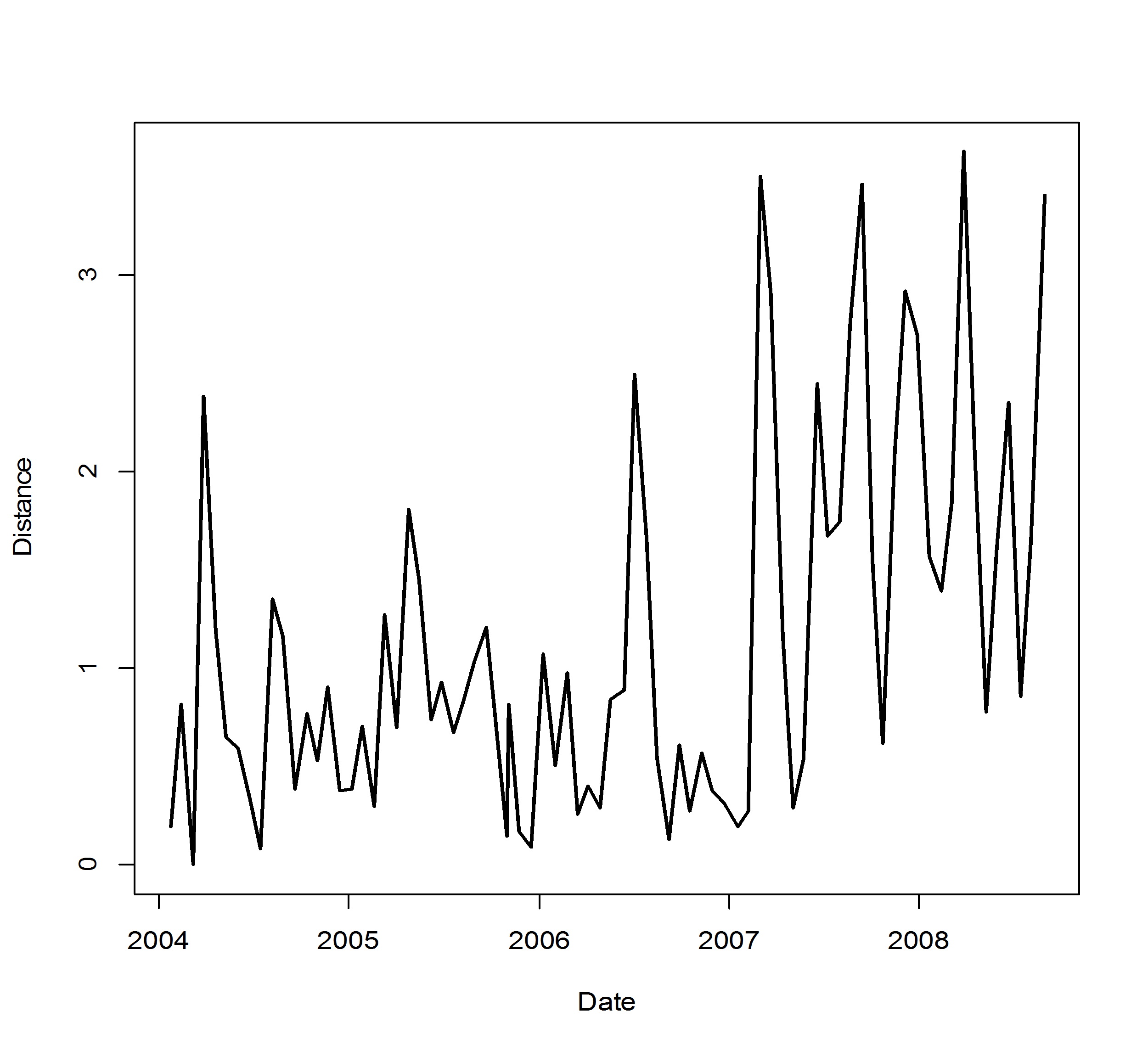}
  \caption{Wasserstein distances between persistence diagrams of the  whole networks derived from the DJIA stocks.}\label{totalwasspaperzn}
\end{figure}

 In \cite{11}, Gidea applied his method on the network derived from the DJIA stocks listed as of February 19, 2008, and the data considered   corresponds to the period between January 2004 and September 2008. The  corresponding graphical representation of the  time series $(X_t)_t$ is presented in Figure \ref{totalwasspaperzn}.  
One can notice from Figure \ref{totalwasspaperzn}, as Gidea concluded in \cite{11},  that the  obtained time series   shows  significant  changes  in  the  topology  of  the  correlation  network  in  the period prior to the onset of the 2007-2008 financial crisis (see \cite{11}).\\

When the network is large, as for instance with social networks, dealing with the entire network lead to a very high computational burden.  It is therefore important to find a way to significantly reduce the execution time of the above method without deteriorating its performance.  To this end, we propose a simplified method  based on the closeness centrality notion.

\section{Persistent homology of dynamic central subnetworks}\label{Sec-Central}

In order to reduce the computation cost of Gidea's method,  one could think of considering a particular subnetwork rather than the whole network.   But, to keep similar performance as with the initial method, the  chosen  subnetwork has to ``summarize'' as well as possible the global network.  Our aim is to  show that the notion of closeness centrality can help  in determining  suitable subnetworks which  we call central subnetworks.\\

In this section, we  start by explaining how the closeness centrality can be used to determine the desired central subnetwork.  We will give an algorithm which details the different steps of extracting the central subnetwork from the whole network. Then,  we show how to  get persistence diagrams of dynamic central subnetworks.\\

 Recall  that   the closeness centrality of a node $i$ in a  weighted network $ G(V, E,\omega)$, with the weighted matrix  $W=(\omega(i,j))$,   is the following quantity:
             $$ C_{c}^{\omega}(i)=\frac{1}{\sum_{j}^{}\delta(i,j)}  $$
             where
             $$\delta(i,j)= min(\frac{1}{\omega(i,h)}+\cdots + \frac{1}{\omega(h,j)}) $$
             such that  $ h $ are intermediary nodes on paths between $i$ and $j$ (see \cite{9} for more details).
   The closeness centrality of a graph helps in detecting nodes that are able to spread information efficiently through the whole graph. Then, based on this fact, we propose that  the central subnetwork would be the one which contains a node having the highest closeness centrality. This node will be called central.\\
   
For practical implementation,  we propose to chose the desired central subgraph from the subgraphs $ G (s) $ defined in  Subsection \ref{PH}. Recall that $ G (s) $ is the subgraph resulting from the sublevel set  $\omega^{-1} ((-\infty,s])$ of the weight function $\omega$  at a threshold  $s$.  Therefore, in order to ensure that the central node is present in the central subgraph, the thresholds $ s $ will be chosen among the weights of the edges incident to the chosen central node. Thus, once a threshold $ s $  is adopted, one can implement the above method, using the matrices  $W=(W_{i,j})$ and 
$\widetilde{W}=(\widetilde{W}_{i,j})$,  as explained at the end of Subsection \ref{implem}. Notice that the list of  weights of the edges incident to the chosen central node form exactly  the line of the weighted matrix   $W$ corresponding to the central node. This row, seen as a  statistical  series, will be sorted in ascending order.  Clearly, its smallest  element is always zero, so it will be not considered and    the resulting ordered  statistical series will be denoted by $ U_c $. Then, the choice of the threshold is determined by the choice of its rank in $ U_c $ which will be called the threshold rank. Therefore, 
 the   algorithm  representing  the steps of constructing the weighted  matrix $\widetilde{W}$ of the  central subnetwork from the matrix $W$ of the  whole network is given in Algorithm 1. 

\begin{algorithm}[H]
\caption{Construct the matrix $\widetilde{W} $}  
\begin{algorithmic}
\REQUIRE Weighted matrix $W=(W_{i,j})$, threshold rank $r$ and  the maximum value of weights $m$.
\ENSURE $ \widetilde{W}=(\widetilde{W}_{i,j})$ The weighted matrix of the central subnetwork
\STATE 1. Compute the closeness centrality of all nodes.
\STATE 2. Determine a central node $c$ of the graph $G$.
\STATE 3. Extract the row   $u_c=(W_{c,j})_{j}$.
\STATE 4. Induce the vector  $U_c=(W_{c,j})_{j\neq c}$.
\STATE 5. Compute the threshold $s$ of rank $r$ in the vector $U_c$.
\STATE 6. $\widetilde{W} \leftarrow W $
\FOR{$0\leq i,j \leq n-1$}
\IF{$W_{i,j}\geq s $}
\STATE $\widetilde{W}_{i,j} \leftarrow m+1$
\ELSE
\STATE $\widetilde{W}_{i,j} \leftarrow W_{i,j} $
\ENDIF
\ENDFOR
\end{algorithmic}
\end{algorithm}
%%%%%%%%%%%%%%%%%%%%%%%%%%%%%%%%%%%%%%%%%%%%%%%%%%%%%%%%%%%%%%%%%%%%%%%%
 
  Now, to obtain the desired time series of $2$-Wasserstein distances between the persistence diagrams of the studied dynamic network, we implement  Algorithm 2.

\begin{algorithm}[H]
\caption{The time series of $2$-Wasserstein distances}
\begin{algorithmic}
\REQUIRE The matrix $\widetilde{W}_t$ obtained using  Algorithm 1 for each time $t$, set  a reference time $t_0$.
\ENSURE   The time series  $(\widetilde{X}_t)_t:=( W_2 (\widetilde{\Delta}_t,\widetilde{\Delta}_{t_0} ) )_t$ between each persistence diagram $\widetilde{\Delta}_t$ and the persistence diagram $\widetilde{\Delta}_{t_0}$ at the time $t_0$
\STATE 1.  Using Rips filtration, create the persistence diagram   $\widetilde{\Delta}_t$ of each time $t$ associated to the  matrix $\widetilde{W}_t$.
\STATE 2. Compute the  $2$-Wasserstein distances $(\widetilde{X}_t)_t$.
\end{algorithmic}
\end{algorithm}
  
\begin{figure}
\centering
\begin{tikzpicture}[node distance = 1.5cm, align= center]
\node (start) [startstop] {Start};
\node (in1) [io, below of=start] {Dynamic network\\ };
\node (pro1) [process, below of=in1] {Algorithm 1};
\node (out1) [io,below of=pro1] {Dynamic central subnework\\ };
\node (pro2) [process, below of=out1] {Algorithm 2};
\node (out2) [io,below of=pro2] {Wasserstein distances on the dynamic central subnework\\ };
\node (stop) [startstop, below of=out2] {Stop};
\draw [arrow] (start) -- (in1);
\draw [arrow] (in1) -- (pro1);
\draw [arrow] (pro1) -- (out1);
\draw [arrow] (out1) -- (pro2);
\draw [arrow] (pro2) -- (out2);
\draw [arrow] (out2) -- (stop);
\end{tikzpicture}
\caption{Flowchart of the proposed method.}  \label{flowchart}
\end{figure}
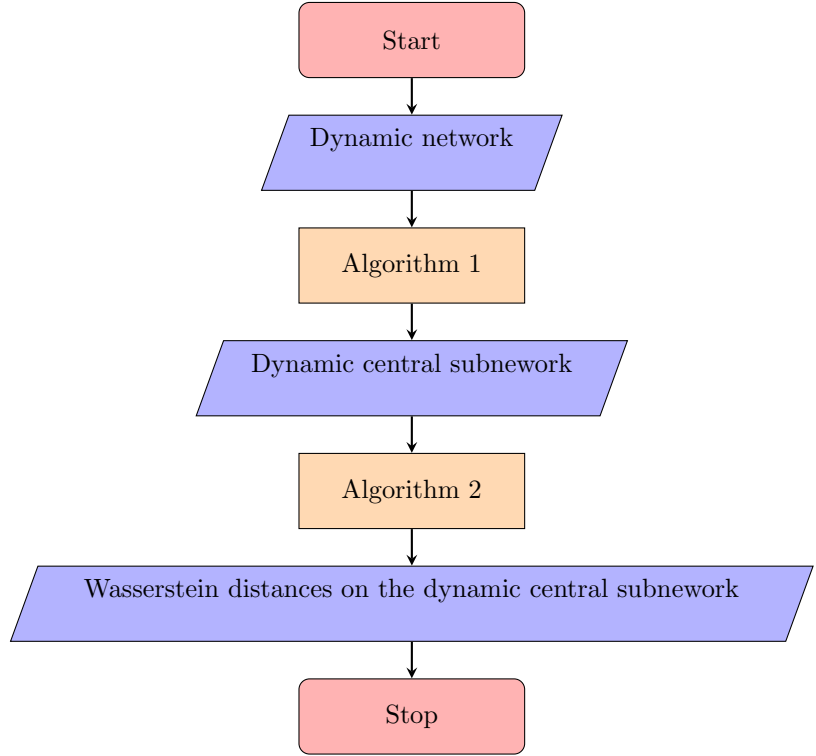

The flowchart of the proposed method is given in Figure \ref{flowchart}.\\

 It is clear that our method is mainly based on the choice of the adequate threshold. To get clear idea about the effectiveness of the choice of threshold that provides good results, we will implement the method for various threshold values. We will focus on five  threshold values corresponding to  remarkable ranks. Explicitly, we exploit the following cases:

     \begin{itemize}
     \item  Min threshold:   The first threshold, $s_0$, will be the   minimum value of the statistical series $ U_c $. 
     \item $Q1$ threshold:  $s_1$ will be the first quartile of $ U_c $.
     \item $Q2$ threshold: $s_2$ will be the second quartile (median) of $ U_c $.
     \item $Q3$ threshold: $ s _3 $ will be the third quartile of the statistical series $ U_c $.
     \item Max threshold: $ s _4  $ will be de the maximum value of the statistical series $ U_c $.\\
   \end{itemize}

\section{Numerical experiments}\label{sec-experimnts}
 In this section, we apply our method described in Section \ref{Sec-Central} to various simulated dynamic networks. \\
 
 We study in each example from which threshold we obtain a good simplification of the whole graphs while maintaining the global observations that we can obtain from the associated  time series. For the latter, in addition to a visual comparison between the  studied time series  charts of $(X_t)_t$ and $(\widetilde{X}_t)_t$, we use the adjusted R-squared. In fact, the reason behind using the adjusted R-squared is based on several observations which indicate that there is some linearity between the studied time series. This hypothesis is confirmed through various examples. 

    %%%%%%%%%%%%%%%%%%%%%%%%%%%%%%%%%%%
%%%%%%%%%%%%%%%%%%%%%%%%%%%%%%%%%%%%%%%%%%%%%%%%%%%%
%%%%%%%%%%%%%%%%%%%%%%%%%%%%%%%%%%%%%%%

  \subsection{First experiment: dynamic network with a fixed central node}\label{sub:experimnt1}
  
For this first experiment, we propose to study a dynamic network with $100$ networks of $51$ nodes and having the same central node.\\

 \noindent \textbf{Description}: As we know, a weighted network is represented by weighted matrix. Then, to simulate a dynamic network with $100$ networks, we simulate $100$ square weighted matrices   of order  $51$. In this experiment we propose to fix the central node in all networks of the studied  dynamic network. To do so, each matrix (representing a weighted  network) will be obtained as follows:
 \begin{itemize}
 \item First, we  simulate $50$ vectors of length $10$ as realizations of a reduced centred normal distribution. 
  \item To isolate the fixed central node during all  $100$ networks, we  take the $51^{st}$  vector as  sum of the other $50$ vectors.
    \item We then compute the correlation matrix $(c_{i,j})$ between these vectors. 
    \item Finally, the desired weighted matrix is  the distance matrix $(d_{i,j})$, where $d_{i,j}=\sqrt{2(1-c_{i,j}}) $. 
 \end{itemize}

Over the 100 networks,  the  $51^{st}$  node  is always the   central node,  it is thus considered to compute the thresholds $s_0$ to $s_4$ as explained at the end of Section \ref{Sec-Central}. 
 
Next,  from each square weighted matrix,  we derive its persistence diagram of
$0$-dimensional features.
 At the end, we compute the times series $(X_t)_t$ of the $2$-Wasserstein distances between all the
obtained $0$-dimensional persistence diagrams and a fixed $0$-dimensional persistence diagram.\\

\begin{figure}[h!]
	\begin{minipage}[b]{0.40\linewidth}
		\centering \includegraphics[scale=0.3]{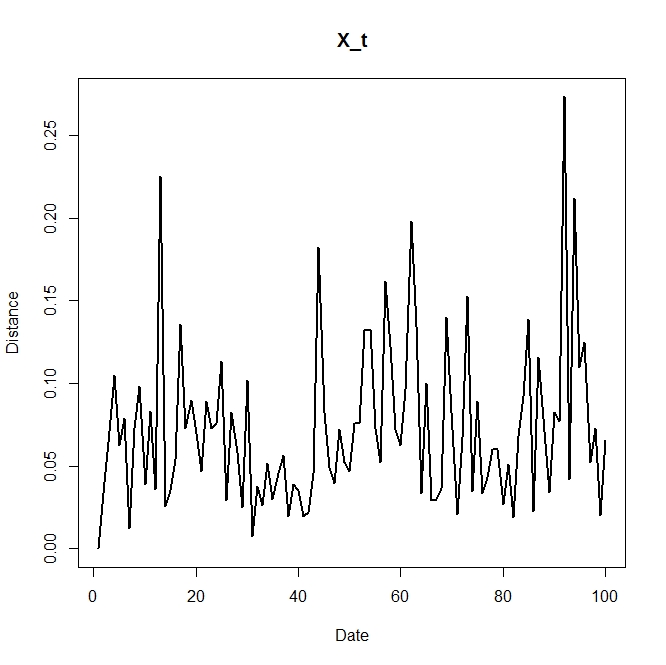}
	\caption{ The time series $({X}_t)_t$   for the first simulation.}\label{fig-ntotalwassim51}
	\end{minipage}\hfill
	\begin{minipage}[b]{0.48\linewidth}	
		\centering \includegraphics[scale=0.3]{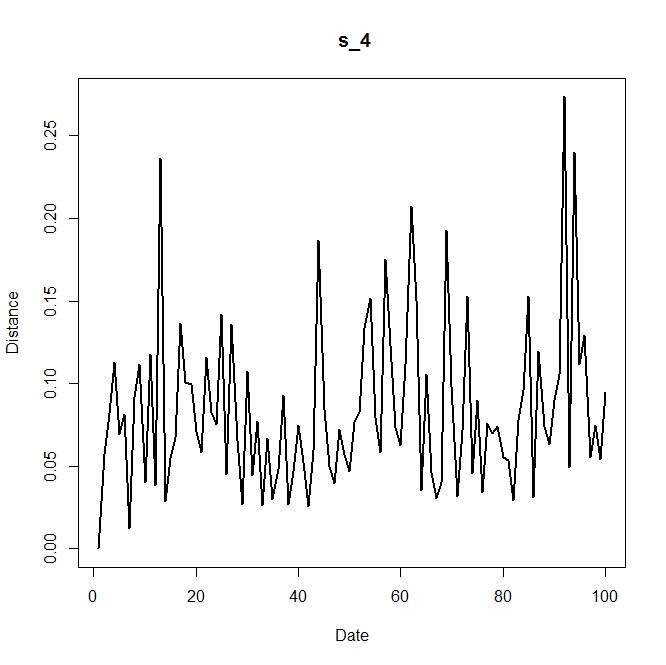}
		 \caption{The time series $(\widetilde{X}_t)_t$ at the threshold $s_4$ for the first simulation.}\label{fig-ncloswassim4q}
	\end{minipage}
	\end{figure}

The graphical representation of this time series  is displayed in Figure \ref{fig-ntotalwassim51}. \\
 
\begin{figure}[h!]
	\begin{minipage}[b]{0.40\linewidth}
		\centering \includegraphics[scale=0.3]{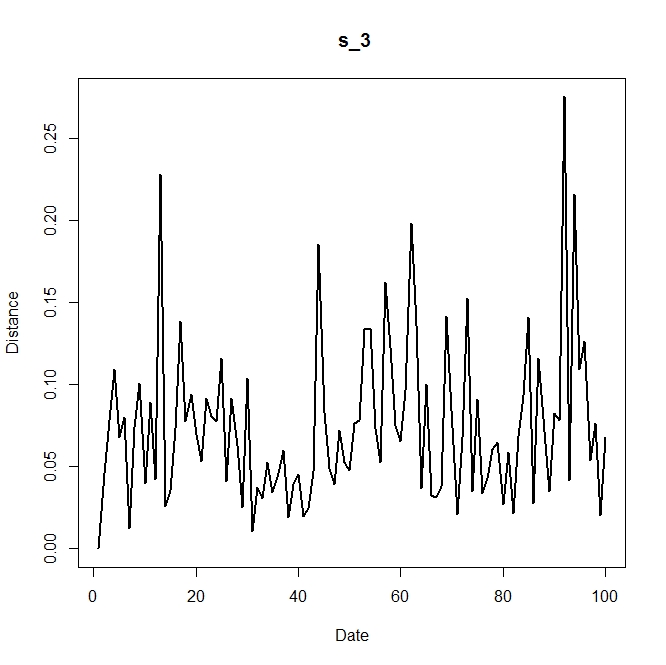}
		 \caption{The time series $(\widetilde{X}_t)_t$ at the threshold $s_3$ for the first simulation.}\label{fig-ncloswassim3q}
	\end{minipage}\hfill
	\begin{minipage}[b]{0.48\linewidth}	
		\centering \includegraphics[scale=0.3]{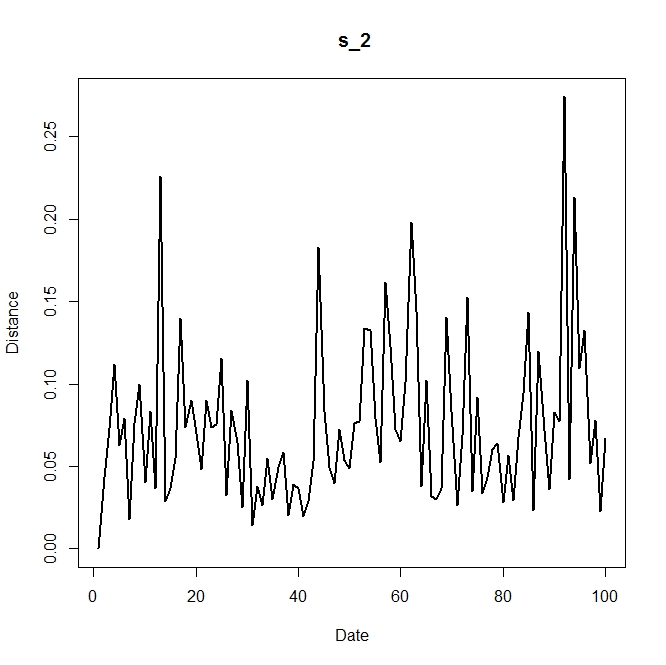}
		 \caption{The time series $(\widetilde{X}_t)_t$ at the threshold $s_2$ for the first simulation}\label{fig-ncloswassim1q}
	\end{minipage}
                    %  \caption{  The time series $({X}_t)_t$  and $(\widetilde{X}_t)_t$ with thresholds $s_3$, $s_2$ and $s_1$  for the DJIA stocks data.}
\end{figure}

\begin{figure}[h!]
	%\begin{minipage}[b]{0.40\linewidth}
		%\centering \includegraphics[scale=0.3]{ntotalwassim51}
	%\caption{ The time series $({X}_t)_t$   for the irst simulation.}
	%\end{minipage}\hfill
	\begin{minipage}[b]{0.48\linewidth}	
		\centering \includegraphics[scale=0.3]{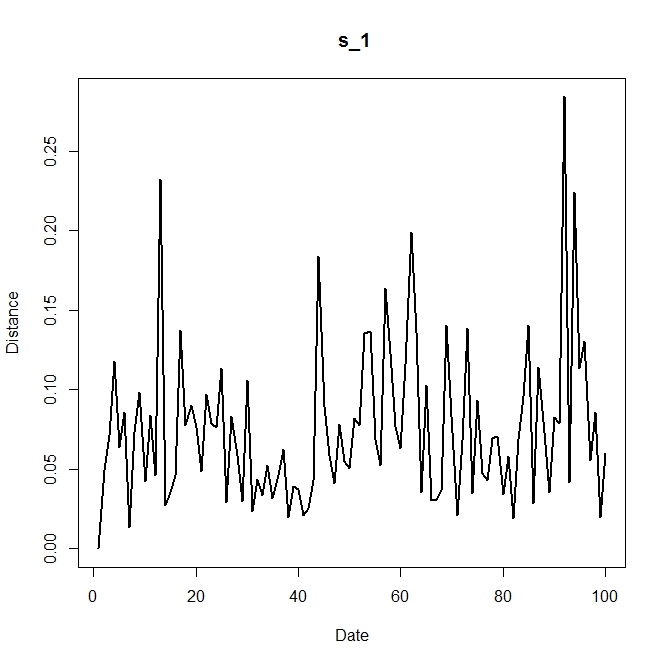}
		 \caption{The time series $(\widetilde{X}_t)_t$ at the threshold $s_1$ for the first simulation.}\label{fig-ncloswassim1q}
	\end{minipage}
	\begin{minipage}[b]{0.40\linewidth}
		\centering \includegraphics[scale=0.3]{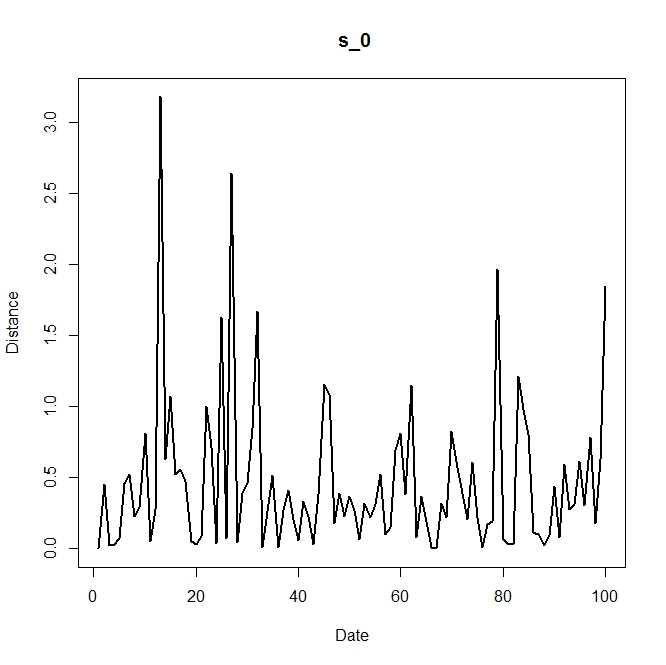}
		 \caption{The time series $(\widetilde{X}_t)_t$ at the threshold $s_0$ for the first simulation.}\label{fig-ncloswassim0q}
	\end{minipage}\hfill
 
\end{figure}

 Now, to see the efficiency  of our approach described  in Section \ref{Sec-Central}, we compute, at   each of the    five   thresholds proposed  at the end of Section  \ref{Sec-Central},  the time series $(\widetilde{X}_t)_t$ associated to the central subnetworks   (see Figures \ref{fig-ntotalwassim51} to \ref{fig-ncloswassim0q}).


%\begin{figure}[H]
  %\centering
  %\includegraphics[scale=0.15]{superpo51tot321n}
%  \caption{Wasserstein distances between persistence diagrams relating to the whole graphs and the central subnetworks thresholds  $s_3$ , $s_2$ and $s_1$for the first experiment.}\label{fig-superpo51tot321n}
%\end{figure}

%\begin{figure}[H]
 % \centering
 % \includegraphics[scale=0.04]{superpo51tot40n}
 % \caption{Wasserstein distances between persistence diagrams relating to the whole graphs and the central subnetworks thresholds  $s_4$  and $s_0$ for the first experiment.}\label{fig-superpo51tot40n}
%\end{figure}
 
\medskip

 \noindent  \textbf{Comment}: It is customary for the comparison between two graphs to first go through a visualization.  We can observe from  Figures \ref{fig-ntotalwassim51} to \ref{fig-ncloswassim0q} that  there is     an almost perfect linear correspondence between the time series $(X_t)_t$ and $(\widetilde{X}_t)_t$  at each   threshold $s_1$, $s_2$, $s_3$ and $s_4$. This is confirmed via the adjusted R-squared values in the third column 3 of Table \ref{table:experience1}.  But, this is not the case at the threshold $s_0$  (see   Figure \ref{fig-ncloswassim0q}).  Namely, the adjusted R-square coefficient at this threshold  is $0.002$.\\

Now, to compare the execution time of each method, we compute the time ratio between them; i.e., the execution time of our algorithm over the execution time of Gidea's one  (see the second column in  Table \ref{table:experience1}).  We can deduce that the execution time is considerably reduced using our algorithm, namely for the thresholds $s_1$, $s_2$ and $s_3$. Notice that the  average percentage of pruned edges, at each threshold, is presented in the fourth column of Table \ref{table:experience1}. Naturally, this average percentage increases when the threshold decreases.

\begin{table}[H]
  \begin{center}
\begin{tabular}{|c|c|c|c|}
  \hline
Threshold    & Time ratio & Adjusted R-squared & Average percentage of pruned edges \\
  \hline
   $s_4$ &0.95 & 0.9389 & 4.17 \\
   $s_3$ & 0.31 & 0.9962 & 39.83\\
   $s_2$ & 0.2 &  0.9977   & 64.42\\
  $s_1$ & 0.16  &   0.991 & 85.32 \\
  \hline
\end{tabular}
\caption{Execution time ratios,  adjusted R-squared and  average percentages of pruned edges  for the first experiment.}\label{table:experience1}
\end{center}
\end{table}

%%%%%%%%%%%%%%%%%%%%%%%%%%%%%%%%%%%%%%%%%%%%%%%%%%%%%%%%%%%%%%%%%%%%%%%%%%  

  \subsection{Second experiment: dynamic network constructed from random covariance matrices}\label{subsec-second-experiment}
In this subsection, we will simulate a dynamic network of 150 networks, where each network  has 60 nodes, by the use of generated covariance matrices.\\ 

 \noindent\textbf{Description:} We   generate $ 60 $ vectors  of length $ 10 $, each one is simulated from a multivariate normal distribution with zero mean and a  covariance matrix generated, as in \cite{6}, by the use of linearly transformed Beta$(\alpha,\alpha)$ distribution on the interval $(-1,1)$\footnote{For this, we use the function \textit{genPositiveDefMat} of the package \textit{clusterGeneration} of the R software.}.  Next, we compute the weighted matrix $(d_{i,j})$ as in Subsection \ref{sub:experimnt1}. This way we end up with a network of $60$ nodes. By repeating this operation $150$ times, each time with a different generated covariance matrix, we built a dynamic network. This is still a naive dynamic network  because the graphs are simulated independently for this experiment.\\
     
\begin{figure}[h!]
	\begin{minipage}[b]{0.40\linewidth}
		\centering \includegraphics[scale=0.3]{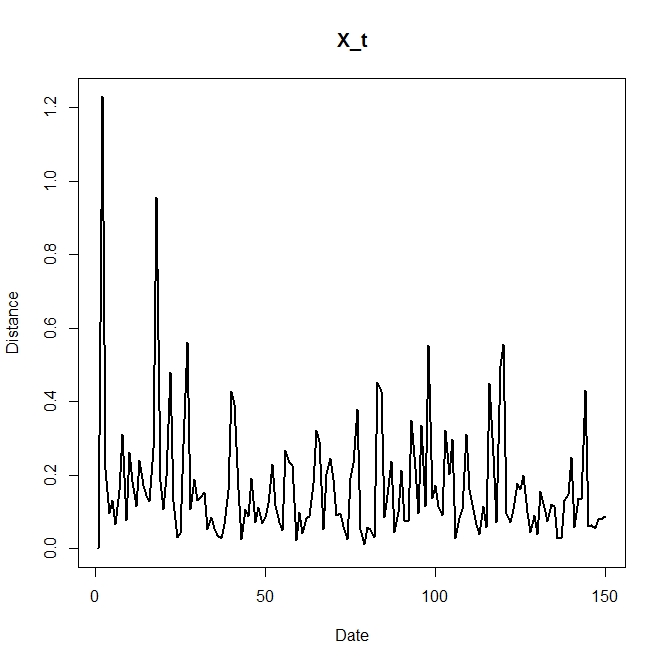}
	\caption{ The time series $({X}_t)_t$   for the  second simulation.}\label{ref-ntotalwassmulti}
	\end{minipage}\hfill
	\begin{minipage}[b]{0.48\linewidth}	
		\centering \includegraphics[scale=0.3]{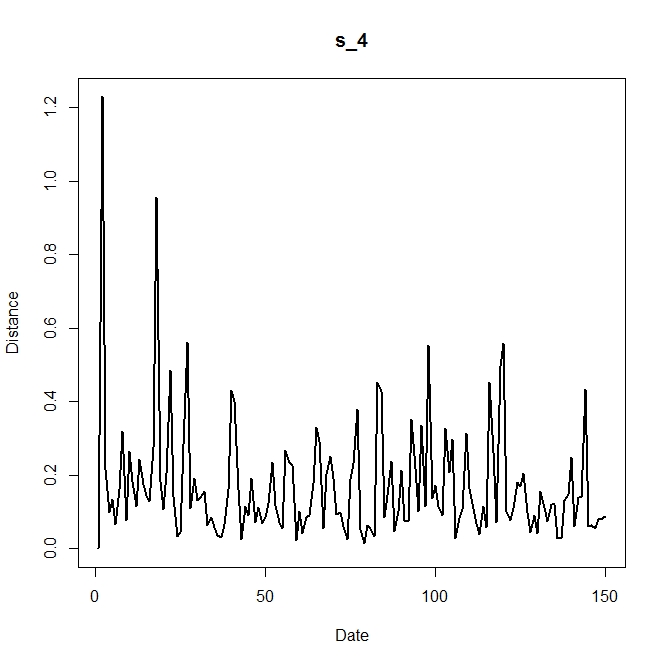}
		 \caption{The time series $(\widetilde{X}_t)_t$ at the threshold $s_4$ for the second simulation.}\label{ref-ncloswassmulti4q}
	\end{minipage}
\end{figure}

The time series chart is displayed in Figure \ref{ref-ntotalwassmulti}. 
 
\begin{figure}[h!]
	\begin{minipage}[b]{0.40\linewidth}
		\centering \includegraphics[scale=0.3]{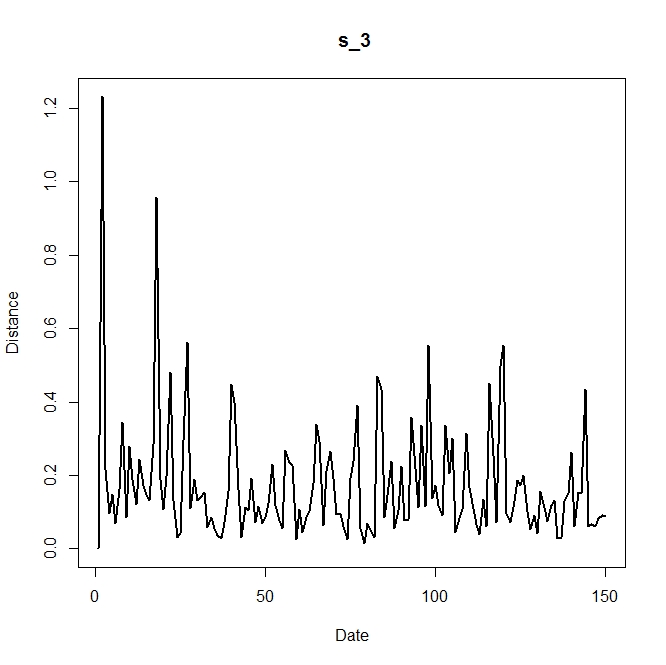}
		 \caption{The time series $(\widetilde{X}_t)_t$ at the  threshold $s_3$ for the second simulation.}\label{ref-ncloswassmulti3q}
	\end{minipage}\hfill
	\begin{minipage}[b]{0.48\linewidth}	
		\centering \includegraphics[scale=0.3]{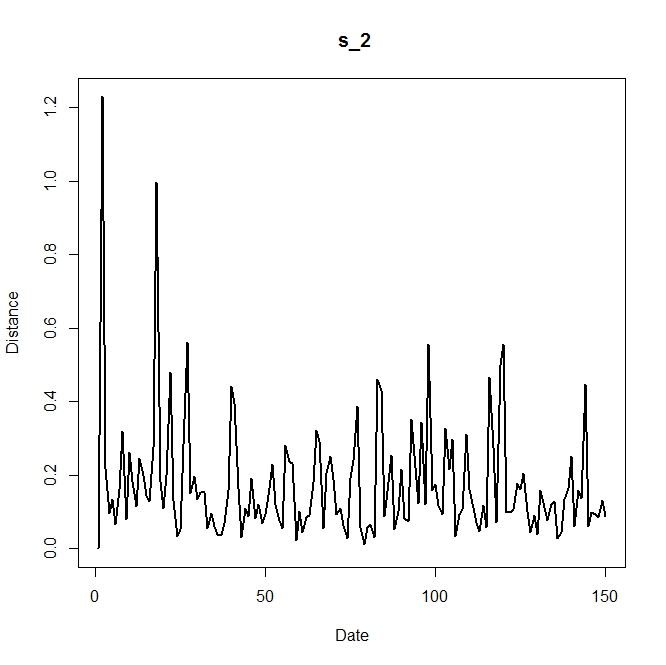}
		 \caption{The time series $(\widetilde{X}_t)_t$ at the threshold $s_2$ for the second simulation.}\label{ref-ncloswassmulti2q}
	\end{minipage} 
\end{figure}

\begin{figure}[h!] 
	\begin{minipage}[b]{0.48\linewidth}	
		\centering \includegraphics[scale=0.3]{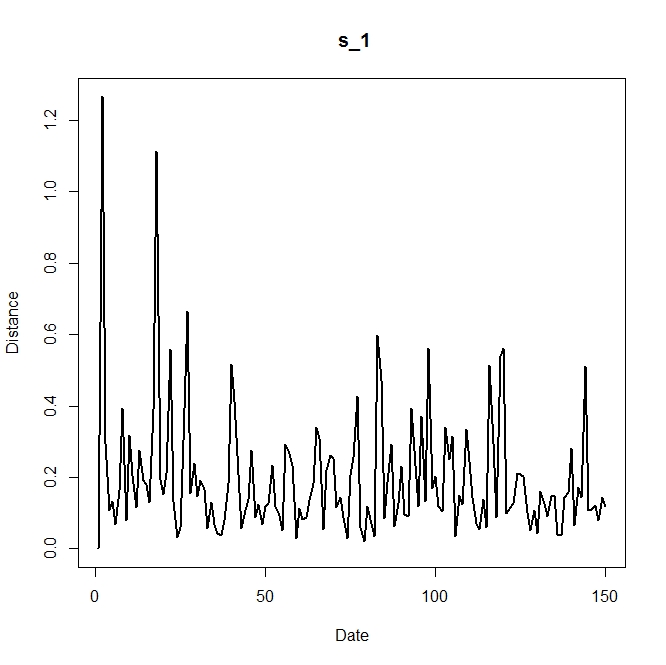}
		 \caption{The time series $(\widetilde{X}_t)_t$ at the threshold  $s_1$ for the second simulation.}
	\end{minipage}
	\begin{minipage}[b]{0.40\linewidth}
		\centering \includegraphics[scale=0.3]{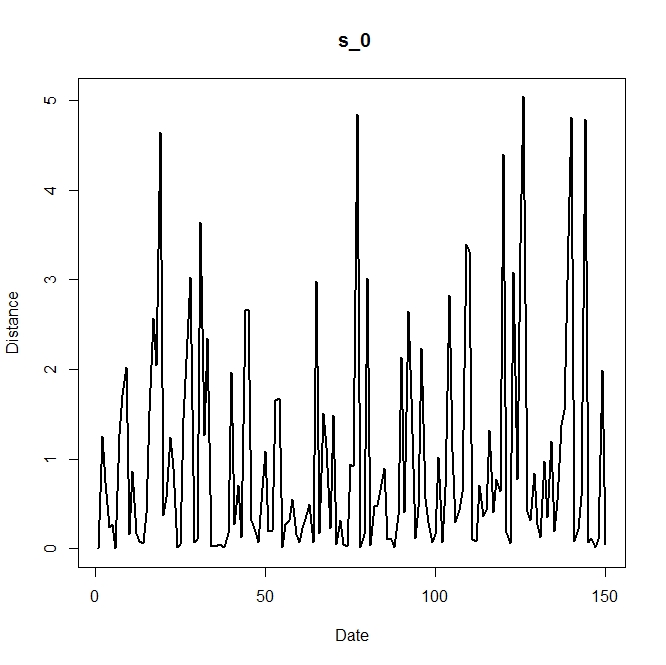}
		 \caption{The time series $(\widetilde{X}_t)_t$ at the threshold $s_0$ for the second simulation.}\label{ncloswassmulti0q}
	\end{minipage}\hfill
\end{figure}

 Now,  we compute, at   each of the    five   thresholds proposed  in Section  \ref{Sec-Central},  the time series $(\widetilde{X}_t)_t$ associated to the central subnetworks   (see Figures \ref{ref-ncloswassmulti4q} to \ref{ncloswassmulti0q}).

As   in Subsection \ref{sub:experimnt1},  we calculate time ratios, adjusted R-squared values and the  average percentages of pruned edges  (see Table \ref{table:experience2}). 

\begin{table}[H]
  \begin{center}
\begin{tabular}{|c|c|c|c|}
  \hline
Threshold    & Time ratio & Adjusted R-squared & Average percentage of pruned edges \\
  \hline
   $s_4$ & 0.86 &   0.9998 & 1.61 \\
   $s_3$ & 0.43 & 0.9987 & 27.83\\
  $s_2$ & 0.22 &  0.997   & 58.85\\
   $s_1$ & 0.16  &  0.9812 & 84.01 \\
  \hline
\end{tabular}
\caption{Execution time ratios,  adjusted R-squared and  average percentages of pruned edges for the second experiment.}\label{table:experience2}
\end{center}
\end{table}

\textbf{Comment:} Although there is no strong central node in the simulated networks for this experiment, we still observe a  good fitting of  $(\widetilde{X}_t)_t$ on $(X_t)_t$, as well a significant improvement  in terms of execution time, in particular with the small   thresholds  $s_1$ and $s_2$.
%%%%%%%%%%%%%%%%%%%%%%%%%%%%%%%%%%%%%%%%%%%%%%%%%%%%%%%%%%%%%%%%%%%%%
  
  %%%%%%%%%%%%%%%%%%%%%%%%%%%%%%%%%%%%%%%%%%%%%%%%%%%%%%%%%%%%%%%%%%%%%%%%%%%%%%%%%

  \subsection{Third experiment: a simulated dynamic network using the central normal law and  an AR(1) process}
In order to test the performance of the method on different sizes of dynamic networks, in this third experiment, we will consider 200 nodes per network. Also, to diversify the simulation methods, we follow a new way to construct the studied dynamic network.  \\
  
\noindent\textbf{Description:}  This  experiment consists of  $200$ weighted matrices, each one of dimension $200 \times 200$, derived from the simulation of 200 vectors of length $20$. The ten first components of each vector are realizations of the reduced centred normal distribution law. The last ten components of each vector come from the AR(1) stationary autoregressive model $ y_t=0.7\times y_{t-1} + \epsilon_t $, where $\epsilon_t$ is a white noise (see \cite{18}). Now,   after calculating  weighted matrices, as in the description of Experiment 1, we end up with a dynamic network.\\

\begin{figure}[h!]
	\begin{minipage}[b]{0.40\linewidth}
		\centering \includegraphics[scale=0.3]{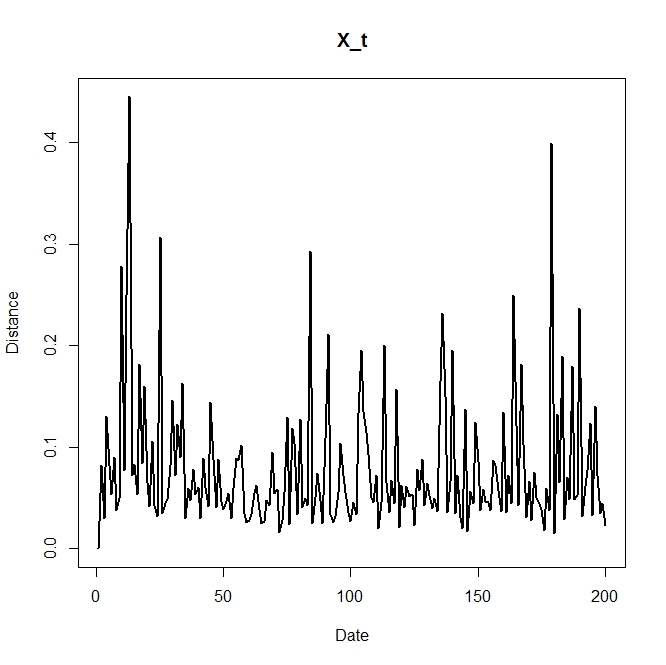}
	\caption{ The time series $({X}_t)_t$   for the  third simulation.}\label{nbigtotalwass}
	\end{minipage}\hfill
	\begin{minipage}[b]{0.48\linewidth}	
		\centering \includegraphics[scale=0.3]{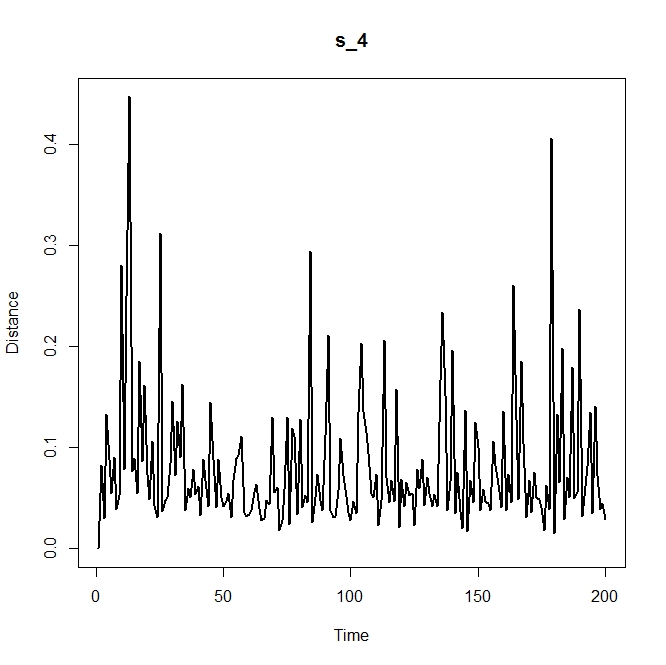}
		 \caption{The time series $(\widetilde{X}_t)_t$ at threshold $s_4$ for the third simulation.}\label{nclosbig4q}
	\end{minipage}
\end{figure}

\begin{figure}[h!]
	\begin{minipage}[b]{0.40\linewidth}
		\centering \includegraphics[scale=0.3]{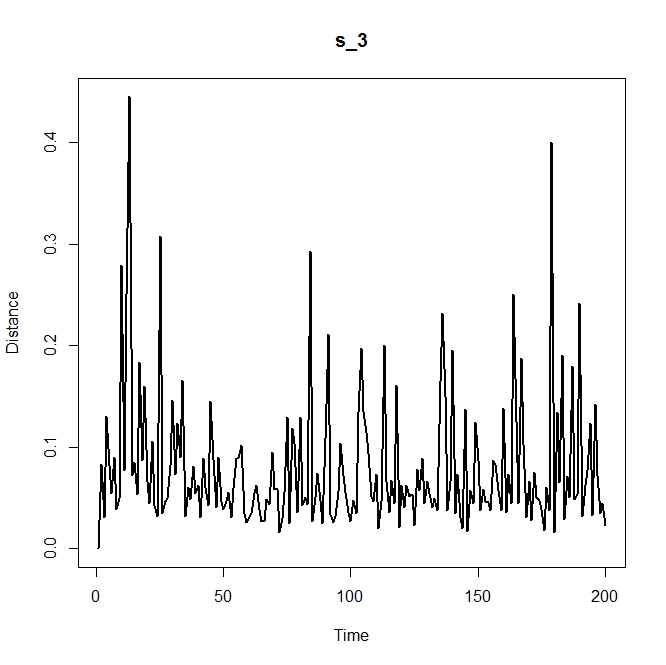}
		 \caption{The time series $(\widetilde{X}_t)_t$ at threshold $s_3$ for the third simulation.}\label{nclosbig3q}
	\end{minipage}\hfill
	\begin{minipage}[b]{0.48\linewidth}	
		\centering \includegraphics[scale=0.3]{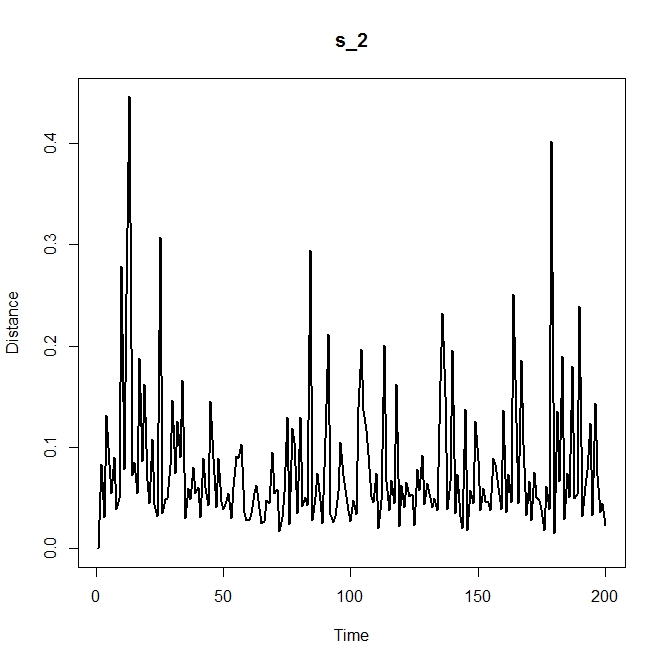}
		 \caption{The time series $(\widetilde{X}_t)_t$ at threshold $s_2$ for the third simulation.}\label{nclosbig2q}
	\end{minipage}
\end{figure}

\begin{figure}[h!]
	\begin{minipage}[b]{0.48\linewidth}	
		\centering \includegraphics[scale=0.3]{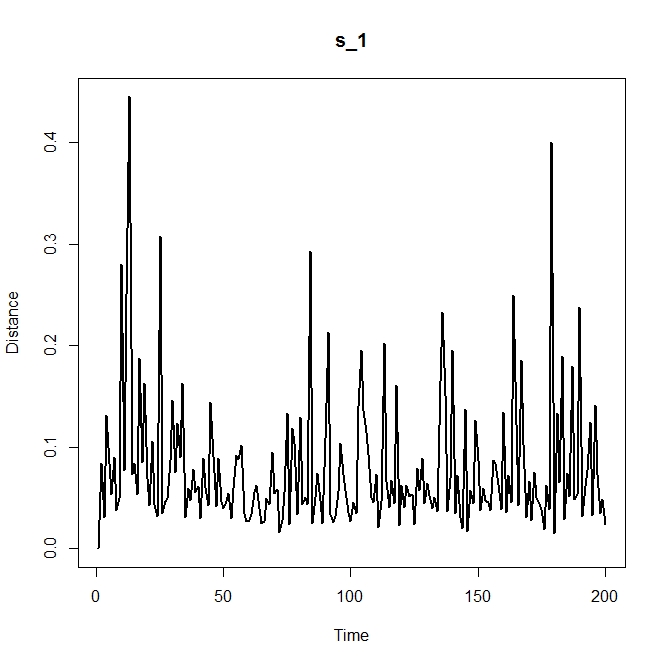}
		 \caption{The time series $(\widetilde{X}_t)_t$ at threshold $s_1$ for the third simulation.}\label{nclosbi1q}
	\end{minipage}
	\begin{minipage}[b]{0.40\linewidth}
		\centering \includegraphics[scale=0.3]{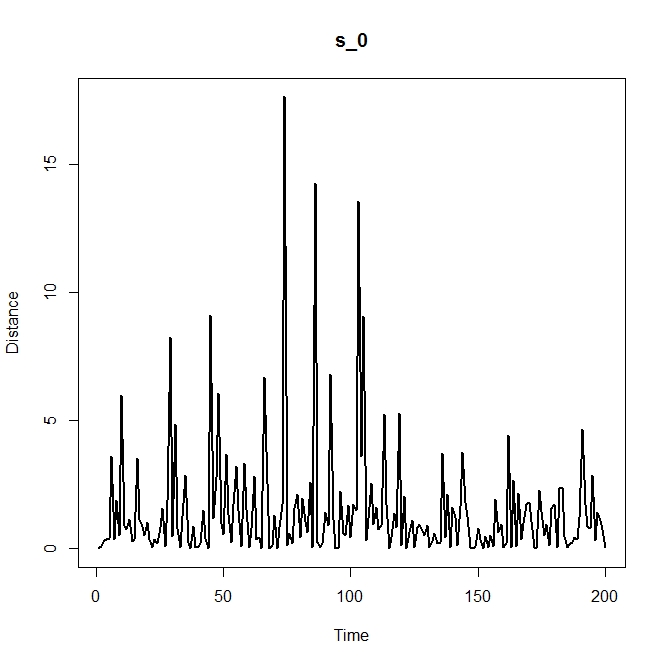}
		 \caption{The time series $(\widetilde{X}_t)_t$ at threshold $s_0$ for the third simulation.}\label{nclosbig0q}
	\end{minipage}\hfill
\end{figure}

The graphical representation of the time series  $(X_t)_t$,  and the one of the central subnetworks $(\widetilde{X}_t)_t$, at the studied thresholds, are  shown  in Figures \ref{nbigtotalwass} to \ref{nclosbig0q}.

Also, we calculate time ratios, adjusted R-squared values and the  average percentages of pruned edges   (see Table \ref{table:experience3}). 

\begin{table}[H]
%\label{table:exp3}
\begin{center}
\begin{tabular}{|c|c|c|}
  \hline
      Threshold    &     Time ratio  &     Adjusted R-squared   \\
  \hline
 $s_4$ & 0.95 & 0.9967  \\
 $s_3$ & 0.3& 0.9998 \\
 $s_2$ & 0.17 &0.9997 \\
$s_1$ & 0.12 &0.9997 \\
  \hline
\end{tabular}
\caption{Execution time ratios,  adjusted R-squared and  average percentages of pruned edges  for the  third experiment.}\label{table:experience3}
\end{center}
\end{table}

\textbf{Comment:} As for Experiment 2, there is no strong  central node in the simulated networks.  We still observe a  good fitting of  $(\widetilde{X}_t)_t$ on $(X_t)_t$ (see Table~\ref{table:experience3}), as well a significant improvement  in terms of execution time, in particular with small thresholds.\\
 Once again,  the information contained in the time series $({X}_t)_t$ is not sufficient when considering the threshold $s_0$ (see Figure \ref{nclosbig0q}).

\subsection{Comparison with edge collapse method}\label{compa-edge}
Recall that the edge collapse algorithm reduces a flag complex to a simpler one with the same persistent homology (see for instance \cite{collapse}). This method must be applied on the $1-$skeleton (graph) of a flag or a clique complex (the same used in our context). So, naturally one could ask for a comparison between this method and ours. This subsection is devoted to the performance comparison between our approach and the edge collapse method, in terms of cost in execution time\footnote{This comparison was suggested by an anonymous referee. We thank her/him for raising this interesting point.}.  For this,  we  apply the edge collapse method as well as our algorithm  on different examples of networks used in this paper.  Namely, we use examples of weighted networks  issued, respectively, from the dynamic networks of Experiments 1, 2 and 3. We use also a weighted network  issued from the dynamic network in Gidea's application (see Subsection \ref{sec-Gidea-appli}) and a new   network with 250 nodes simulated in the same way as the networks in the second experiment. \\

The results are displayed in Table \ref {tbl:compa-e-collaps}.\\

\begin{table}[htbp]
  \centering
  \begin{tabular}{|c|c|c|c|c|c|}
    \hline
    \multicolumn{1}{|c|}{} &\multicolumn{5}{|c|}{Execution time in seconds}  \\ \hline
     \multicolumn{1}{|c|}{Number of nodes}  & $s_4$ & $s_3$ &$s_ 2$ & $s_1$ & $\text{edge collapse}$   \\ \hline
  \multirow{1}{*} 
    28 (Gidea's application)& 0.041 &  0.046 &  0.04 &  0.038 &  0.031 \\\cline{1-6}
    51 (Experience 1) & \(0.16\) &  \(0.08\) &  \(0.05\) &  0.04 &  0.06 \\\cline{1-6}
    60 (Experience 2) & \(0.072\)&   \(0.07\)&   \(0.061\)    & \(0.0575\) &    \(0.0625\) \\ \cline{1-6}
   200 (Experience 3)& \(13.09\) & \(3.73\)& \(1.15\)&   \(0.25\)&   \(95.57\)\\ \cline{1-6}
   250 (New simulation for Experience 2)&  \(31.20\) & \(10.82\)& \(6.10\)&   \(0.47\)&   \(368.5\)\\

   \hline  
  \end{tabular}
  \caption{ Execution times of our method and of the edge collapse method.}
  \label{tbl:compa-e-collaps}
\end{table}

\noindent\textbf{Comment:}  
From Table \ref{tbl:compa-e-collaps}, we see that,  for Experience 1 (51 nodes) and Experience 2 (60 nodes), the edge collapse method is faster than ours at the thresholds $s_4$ and $s_3$ while it is not the case  at the  thresholds $s_2$ and $s_1$. Moreover, our method shows a better execution time for the 200 node graph. \\

Following these observations, it appears that the edge collapse method is more efficient as soon as the number of nodes is small. This is confirmed by the data in the first row of Table \ref{tbl:compa-e-collaps} which represents the results of the applications of the two methods on a graph with 28 nodes resulting from the networks used in  Subsection \ref{sec-Gidea-appli}.  \\

However, as soon as the number of edges is larger it appears that our method is faster.  To confirm this observation, we simulate a network with 250 nodes as in the second experiment (see Subsection \ref{subsec-second-experiment}). 
Therefore, the results of the application presented in the last row of Table \ref{tbl:compa-e-collaps} confirm the assumption. 

Finally,  besides comparing the two methods in terms of execution time, it should be noted that   the edge collapse method provides theoretical guarantees and allows to obtain, in output, the same starting persistence diagram.

\section{Application on dynamic real   data networks}\label{appli dynamic}
In this section, we evaluate our method on two real data networks.

\subsection{Dynamic network derived from the DJIA stocks}\label{appli DJIA}
In Subsection \ref{sec-Gidea-appli},  we recalled how Gidea used topological data analysis to study a  financial network. Here, we show how our method performs faster than Gidea's one  without compromising quality and results. 

\begin{figure}[h!]
	\begin{minipage}[b]{0.40\linewidth}
		\centering \includegraphics[scale=0.3]{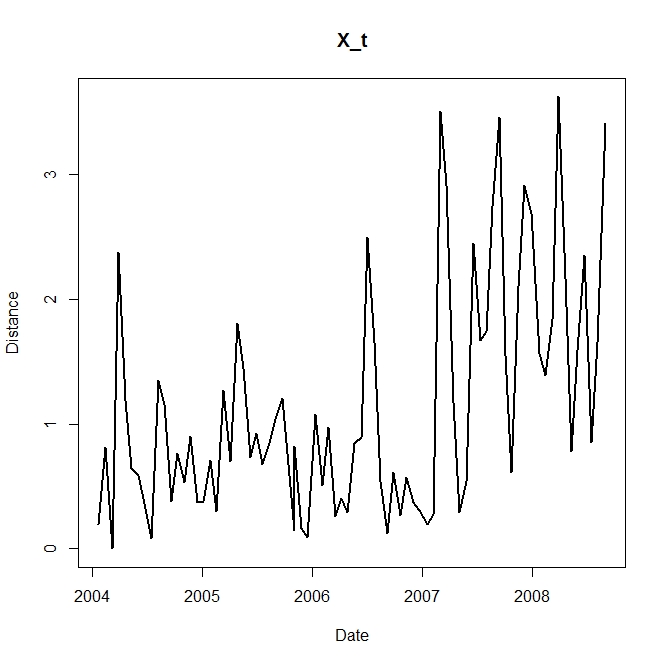}
	\caption{The time series $({X}_t)_t$   for the DJIA stocks data.}\label{ntotalwasspaperz}
	\end{minipage}\hfill
	\begin{minipage}[b]{0.48\linewidth}	
		\centering \includegraphics[scale=0.3]{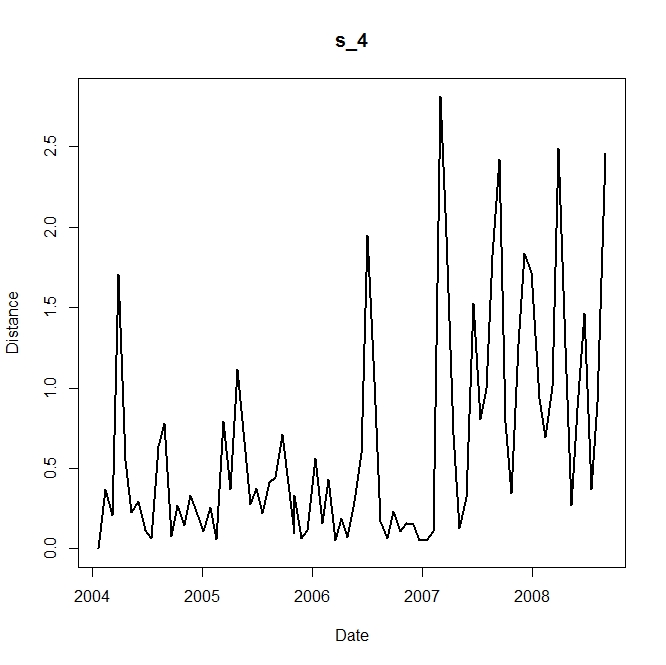}
		 \caption{The time series $(\widetilde{X}_t)_t$ at threshold $s_4$ for the DJIA stocks data.}\label{ncloswasspaper4q}
	\end{minipage}
\end{figure}

\begin{figure}[h!]
	\begin{minipage}[b]{0.40\linewidth}
		\centering \includegraphics[scale=0.3]{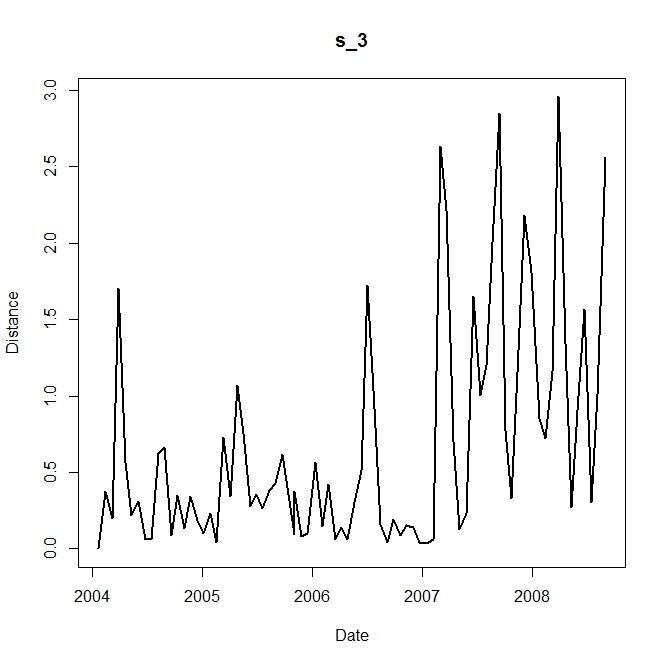}
		 \caption{The time series $(\widetilde{X}_t)_t$ at threshold $s_3$ for the DJIA stocks data.}\label{ncloswasspaper3q}
	\end{minipage}\hfill
	\begin{minipage}[b]{0.48\linewidth}	
		\centering \includegraphics[scale=0.3]{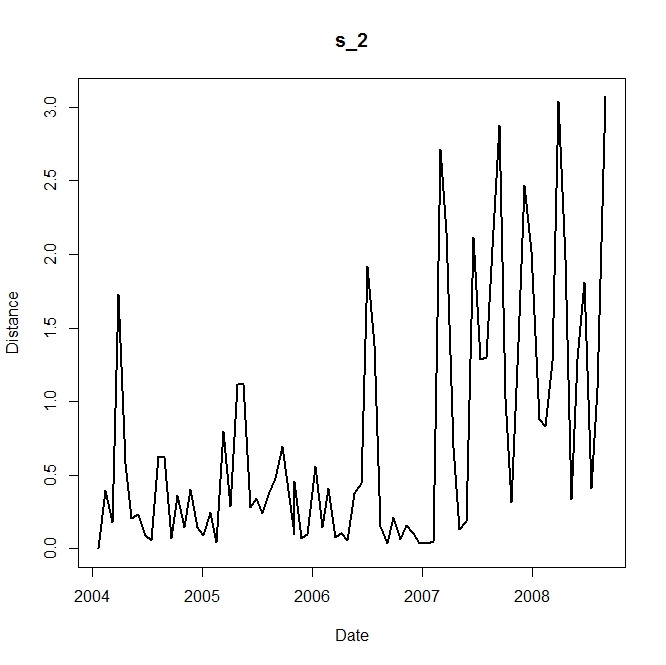}
		 \caption{The time series $(\widetilde{X}_t)_t$ at threshold $s_2$ for the DJIA stocks data.}\label{ncloswasspaper2q}
	\end{minipage}
                    %  \caption{  The time series $({X}_t)_t$  and $(\widetilde{X}_t)_t$ with thresholds $s_3$, $s_2$ and $s_1$  for the DJIA stocks data.}
\end{figure}

\begin{figure}[h!] 
	\begin{minipage}[b]{0.48\linewidth}	
		\centering \includegraphics[scale=0.3]{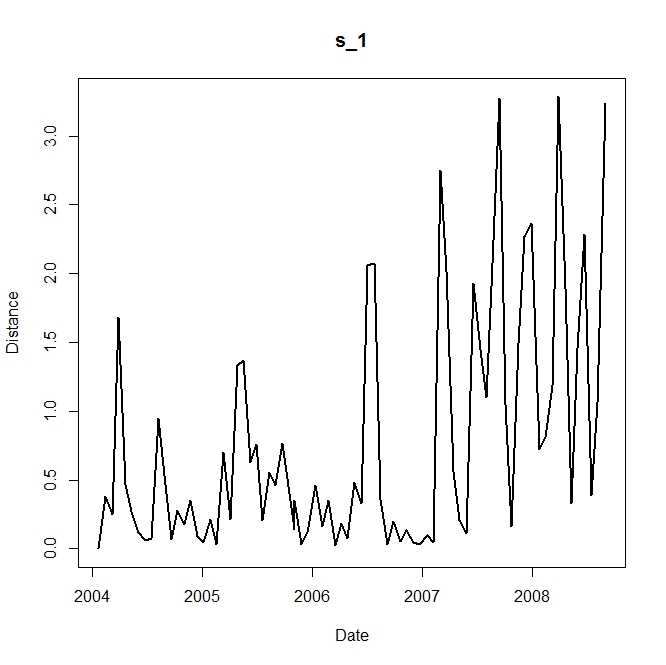}
		 \caption{The time series $(\widetilde{X}_t)_t$ at threshold $s_1$ for the DJIA stocks data.}\label{ncloswasspaper1q}
	\end{minipage}
	\begin{minipage}[b]{0.40\linewidth}
		\centering \includegraphics[scale=0.3]{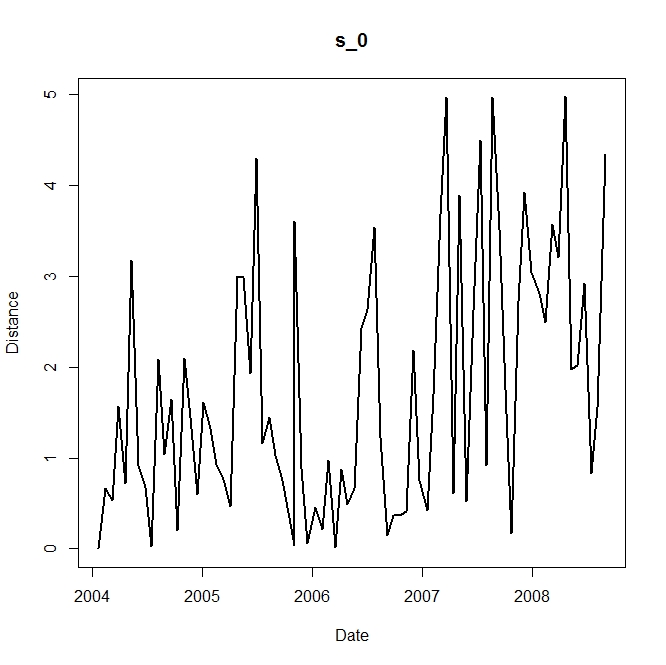}
		 \caption{The time series $(\widetilde{X}_t)_t$ at threshold $s_0$ for the DJIA stocks data.}\label{ncloswasspaper0q}
	\end{minipage}\hfill
\end{figure}

Thus, we apply our method, at the proposed thresholds,   on  the network derived from the DJIA stocks described in Subsection  \ref{sec-Gidea-appli}.   The  charts of these time series are given in Figures \ref{ncloswasspaper4q} to \ref{ncloswasspaper0q}.  The chart of the time series $({X}_t)_t$ for the DJIA stocks data given in Figure \ref{totalwasspaperzn} is put next to the new ones (in Figure \ref{ntotalwasspaperz}) for easy visual comparison.
 
 We can see that  the early signs of a financial market crisis in the time series $(X_t)_t$ corresponding to the whole network is retained at the same period in the new time series $(\widetilde{X}_t)_t$   at the thresholds $s_1$, $s_2$, $s_3$ and $s_4$.  In fact, at these thresholds, one can observe that the time series  $(\widetilde{X}_t)_t$ and $(X_t)_t$ behave  similarly,  as confirmed by the  adjusted R-squared value (see Table \ref{table:rspaper}), although we have pruned a significant number of edges from the whole network (see the third column of Table \ref{table:rspaper}).

\begin{table}[H]
  \begin{center}
\begin{tabular}{|c|c|c|}
  \hline
 Threshold    & Adjusted R-squared & Average percentage of pruned edges  \\
  \hline
   $s_4$ &  0.9933 & 8.55 \\
   $s_3$&   0.9946 & 56.62 \\
   $s_2$ &      0.9765 & 78.87 \\
  $s_1$ &      0.9371 &  92.14 \\
  \hline
\end{tabular}
\caption{Adjusted R-squared in the linear regression of the times series $(X_t)_t$ on $(\widetilde{X}_t)_t$ for the DJIA stocks data.}\label{table:rspaper}
\end{center}
 
\end{table}

Table \ref{stock} gives the time ratio corresponding to the proportion of the execution time of our approach over the execution time of   Gidea's one.  Our approach reduced the execution time by $21\%$, $61\%$, $72\%$ and $74\%$ at thresholds $s_4$, $s_3$, $s_2$ and $s_1$ respectively.

\begin{table}[H]
  \begin{center}
\begin{tabular}{|c|c|}
  \hline
Threshold    & Time ratio  \\
  \hline
  $s_4$ & 0.79   \\
  $s_3$ & 0.39   \\
   $s_2$ & 0.28    \\
   $s_1$ &   0.26      \\
  \hline
\end{tabular}
  \caption{Execution time ratios for the application on DJIA stocks data. }\label{stock}
\end{center}
\end{table}
However, when we consider the threshold $s_0$ (see the graph in  Figure \ref{ncloswasspaper0q}),   one can notice   that   the time series  $(\widetilde{X}_t)_t$ and $(X_t)_t$ are clearly different, as confirmed by the value of the adjusted R-squared in the linear regression of the times series $(\widetilde{X}_t)_t$ on $(X_t)_t$ which is $0.45$.  Nevertheless, the slight deviation in the behaviour of the time series on the $s_0$ scale  could reveal topological properties of the network at the local level which are not detected by the global network.
\\

In conclusion, based in the above observations, the optimal  threshold of our method lies between the first  and the third quadrant of the statistical series of weights of the edges incident to the central node.

\subsection{Cryptocurrency dynamic  network}\label{crypto}
Cryptocurrency dynamic  networks have been a subject of studies of several works (see for instance \cite{12,25,3}).  Here, we are interested in  a dynamic network constructed from the multivariate time series of the closing prices of the four cryptocurrencies Bitcoin, Ethereum, Litecoin and Ripple, from August 24, 2016 to February 19, 2020\footnote{The data was   downloaded from the webpage \url{www.investing.com}}. The dynamic network is constructed, following  \cite{25}, as follows: Let $c_{i}(t)$	 be the closing price of the $i$-th cryptocurrency at the trading day $t$. The Log-return $r_{i,t}$ at the trading day $t$ is defined as $ln(\frac{c_{i}(t)}{c_{i}(t-1)})$ for $i=1,2,3,4$ and $t=1,...,1274$. 
Each trading day $t$ is mapped to the point $x_t=(r_{1,t},r_{2,t},r_{3,t},r_{4,t})\in \mathbb{R}^4$. We end up with the point cloud $ X=\{x_t\}_{t=1}^{t=1274}$ embedded in
 $\mathbb{R}^4$. Now, applying a sliding window of size $w=50$ to the point cloud $X$, we obtain 1225 time varying point clouds $X(t)=\{x_{t-49},...,x_t\}_{t=50}^{t=1274}$. Therefore, the desired weighted dynamic network is $(X(t))_t$ endowed with the weight function defined in Subsection \ref{sec-Gidea-appli}.\\

\begin{center}
  \begin{figure}[H]
 \centering
  \includegraphics[scale=1]{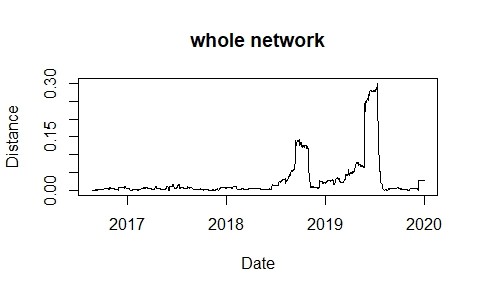}
  \caption{Time series $ ({X} _t) _t $  of  the four cryptocurrencies  stocks data.}\label{fig-wholecc}
\end{figure}
\end{center}

 Now for each point cloud $X(t)$, we  compute the persistent homology of its Rips filtration and  derive the corresponding persistence diagram. 
Subsequently, the time series of the 2-Wasserstein distances between all the persistence diagrams and a fixed one are calculated. The   time series chart is given in  Figure \ref{fig-wholecc}. 

This time series    shows  two periods of behavioural change. The first period was marked by a  first peak recorded on September 14, 2018, followed by a near-recession before declining again on November 1, 2018. The second period is where the time series is increased between May 27, 2019 and July 13, 2019. Notice that there are two major incidents in the cryptocurrency financial market. The first one happened on November 07, 2018 and the second one on June 26, 2019 (see \cite{25}). The first peak of each period can be considered as an early sign of a crisis \cite{11}. \\

\begin{figure}[h!]
	\begin{minipage}[b]{0.40\linewidth}
		\centering \includegraphics[scale=0.3]{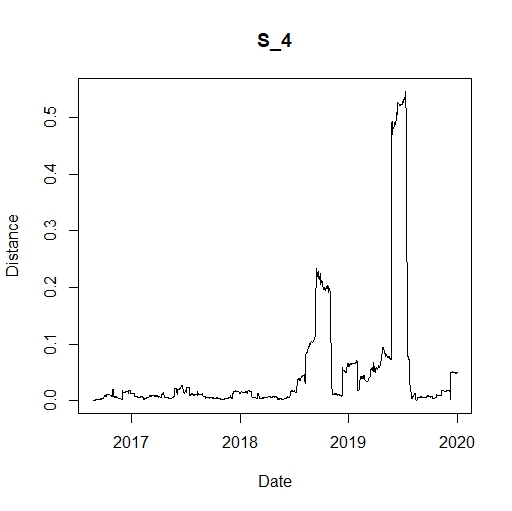}
		%\caption{ Légende 1}
	\end{minipage}\hfill
	\begin{minipage}[b]{0.48\linewidth}	
		\centering \includegraphics[scale=0.3]{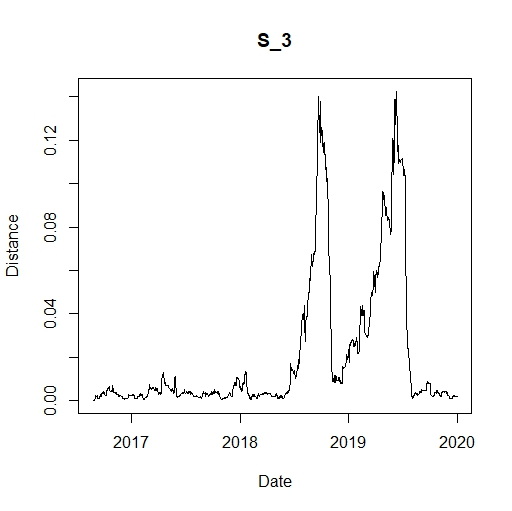}
		%\caption{Légende 2}
	\end{minipage}
	\begin{minipage}[b]{0.40\linewidth}
		\centering \includegraphics[scale=0.3]{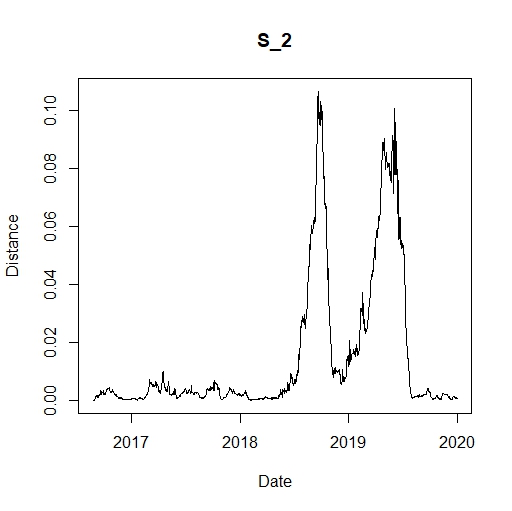}
		%\caption{Légende 3}
	\end{minipage}\hfill
	\begin{minipage}[b]{0.48\linewidth}	
		\centering \includegraphics[scale=0.3]{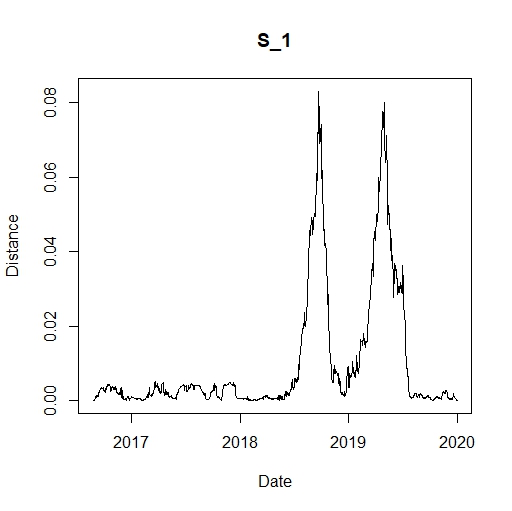}
		%\caption{Légende 4}
	\end{minipage}
                         \caption{Time series $ (\widetilde {X} _t) _t $  concerning the four cryptocurrency  stocks data, of central  subnetworks at, respectively, the thresholds $s_4$, $s_3$, $s_2$ and $s_1$.}\label{fig-central-lecc}
\end{figure}

Now, we apply our algorithm on this dynamic network. The time series charts are given in  Figure \ref{fig-central-lecc}. \\

As in the previous examples,  there is a strong linear relationship between the time series  $(X_t)_t$ and $(\widetilde{X}_t)_t$, especially at thresholds $s_4$ and $s_3$. 
This is confirmed by the adjusted R-squared values (see the first column of Table \ref{crypto}). Moreover, one can observe that there is a positive correlation between the two time series. To confirm this observation, we use  Kendall's tau. Indeed,  Kendall's tau makes it possible to measure  the strength and direction of association that exists between two variables (see \cite{kendall}).  The values of Kendall's tau, listed in the third column of Table \ref{crypto},   allow us to confirm the existence of a strong positive correlation between the two series $(X_t)_t$ and $(\widetilde X_t)_t$, especially at thresholds $s_4$ and $s_3$.\\

\begin{table}[H]
  \begin{center}

\begin{tabular}{|c|c|c|c|}
  \hline
 Threshold    & Time ratio & Adjusted R-squared  & Kendall's tau  \\
  \hline
   $s_4$ & 0.98  & 0.98 &0.80\\
   $s_3$&   0.39& 0.76& 0.64\\
   $s_2$ &      0.23 & 0.55 &0.54\\
  $s_1$ &      0.17  &  0.35 &0.49\\
  \hline
\end{tabular}
  \caption{Execution time ratios,  adjusted R-squared and  Kendall's tau values with the studied thresholds for the cryptocurrencies stock data. }\label{crypto}
\end{center}
\end{table}

Now, to see the efficiency of our algorithm, we compare the execution times via  the execution time ratio as in the previous examples (see the second column  of Table \ref{crypto}).  The  reduction in execution time was  $2\%$, $61\%$, $77\%$ and $83\%$ at thresholds $s_4$, $s_3$, $s_2$ and $s_1$ respectively. \\

In conclusion, based in the above observations, the optimal  threshold of our method lies between the second  and the third quadrant of the statistical series of weights of the edges incident to the central node.

%%%%%%%%%%%%%%%%%%%%%%%%%%%%%%%%%%%%%%
\bigskip

 \section*{Acknowledgements} 
We thank the anonymous reviewers for their careful reading and many insightful comments and suggestions that have helped us to significantly improve the quality of our article.

%%%%%%%%%%%%%%%%%%%%%%%%%%%%%%%%%%%%%

\bigskip

 \section*{Conclusion}  In many context,  weighted networks are very large and their study requires special treatment.
 The proposed method is centred on the simplification of the graphs by using a particular thresholding which is based on the notion of centrality of a graph.\\
 According to the established various simulations as well as  applications on real data, it can be seen that the proposed method has proved its effectiveness for  the threshold $s_3$, and some times even for the thresholds $s_1$ and $s_2$. Thus, these thresholds make it possible to simplify a graph into a central subnetwork which represents it and retains its characteristics. Moreover, several given adjusted R-squared values assure that there is a linear relationship  between the times series $(X_t)_t$ and $(\widetilde X_t)_t$. However, it would be very important to investigate this conclusion from theoretical point of view. Also, the proposed method is based on determining the ``good" threshold which provides considerable simplification. In our study, we discussed ``good"  thresholds related to each example, but determining theoretically the optimal threshold of a given network is an open question.\\
In addition to the open questions above, we present some future research directions  suggested by anonymous reviewers. First, our approach is probably related to persistent local homology \cite{fasy2016exploring}. Indeed, by eliminating some edges and extracting subnetworks from the whole graph in this way, and then computing persistent homology for these simplified graphs; in a sense we are computing ``localized" topological features. Moreover, persistent homology is known to be challenging to compute \cite{fasy2016exploring}. We plan to investigate the interest of our method for computing more efficiently local topological features.  \\
Finally,  distributed persistent homology computation would help also in our purpose. Indeed, several distributed approaches have been proposed to compute persistent homology (see for instance \cite{BKR}  and \cite{MVsW}). We consider that our approach by simplifying a filtration is different from those approaches which organize computation over several clusters. However, these two family of approaches are complementary: after simplifying a filtration, it is also possible to take advantage of distributed algorithms. This also will be considered in our future research works.

\end{document}